\documentclass[11pt]{article}
\usepackage{latexsym,amssymb}
\usepackage{amsmath}
\usepackage[abbrev]{amsrefs}
\usepackage{color}

\newtheorem{Theorem}{\bf Theorem}[section]
\newtheorem{Lemma}{\bf Lemma}[section]
\newtheorem{Proposition}{\bf Proposition}[section]
\newtheorem{Corollary}{\bf Corollary}[section]
\newtheorem{Remark}{\bf Remark}[section]
\newtheorem{Example}{\bf Example}[section]
\newtheorem{Definition}{\bf Definition}[section]
\newtheorem{Comment}{\bf Comment}[section]

\newenvironment{theorem}{\begin{Theorem}$\!\!\!$}{\end{Theorem}}
\newenvironment{lemma}{\begin{Lemma}$\!\!\!$}{\end{Lemma}}

\newenvironment{corollary}{\begin{Corollary}$\!\!\!$}{\end{Corollary}}
\newenvironment{remark}{\begin{Remark}$\!\!\!$}{\end{Remark}}

\newenvironment{definition}{\begin{Definition}$\!\!\!$}{\end{Definition}}

\def\Xint#1{\mathchoice
{\XXint\displaystyle\textstyle{#1}}%
{\XXint\textstyle\scriptstyle{#1}}%
{\XXint\scriptstyle\scriptscriptstyle{#1}}%
{\XXint\scriptscriptstyle\scriptscriptstyle{#1}}%
\!\int}
\def\XXint#1#2#3{{\setbox0=\hbox{$#1{#2#3}{\int}$}
\vcenter{\hbox{$#2#3$}}\kern-.5\wd0}}

\def\dashint{\Xint-}

\numberwithin{equation}{section}

\begin{document}

\title{Existence of solutions to nonlinear parabolic equations\\ via majorant integral kernel}
\author{
\qquad\\
Kazuhiro Ishige, Tatsuki Kawakami, and Shinya Okabe
}
\date{}
\maketitle
\begin{abstract}
We establish the existence of solutions to the Cauchy problem 
for a large class of nonlinear parabolic equations 
including fractional semilinear parabolic equations, 
higher-order semilinear parabolic equations, and viscous Hamilton-Jacobi equations 
by using the majorant kernel introduced in 
[K. Ishige, T. Kawakami, and S. Okabe, Ann. Inst. H. Poincar\'e Anal. Non Lin\'eaire 37 (2020), 1185--1209].  
\end{abstract}
\vspace{25pt}
\noindent Addresses:

\smallskip
\noindent K.~I.:  Graduate School of Mathematical Sciences, The University of Tokyo,\\
\qquad\,\,\, 3-8-1 Komaba, Meguro-ku, Tokyo 153-8914, Japan.\\
\noindent 
E-mail: {\tt ishige@ms.u-tokyo.ac.jp}\\

\smallskip
\noindent 
T. K.: Applied Mathematics and Informatics Course,
Faculty of Advanced Science\\ 
\qquad\,\,\,\,\, and Technology, Ryukoku University,
1-5 Yokotani, Seta Oe-cho, Otsu, \\
\qquad\,\,\,\,\, Shiga 520-2194, Japan.\\
\noindent 
E-mail: {\tt kawakami@math.ryukoku.ac.jp}\\

\smallskip
\noindent S. O.:  Mathematical Institute, Tohoku University,
Aoba, Sendai 980-8578, Japan.\\
\noindent 
E-mail: {\tt shinya.okabe@tohoku.ac.jp}\\
\vspace{20pt}

\newpage
\section{Introduction}
Let $\ell$, $m\in\{0,1,\dots\}$. 
Consider the nonlinear integral equation
\begin{equation}
\tag{I}
\begin{split}
 & u(x,t)=\int_{{\mathbb R}^N}G(x,y,t)\phi(y)\,dy\\
 & \qquad
+\sum_{|\alpha|=\ell}a_\alpha\int_0^t\int_{{\mathbb R}^N}\partial_x^\alpha G(x,y,t-s)
F\left(y,s,u(y,s),\dots,\nabla^m u(y,s)\right)\,dy\,ds
\end{split}
\end{equation}
for $ x\in{\mathbb R}^N$ and $t>0$, 
where 
$N\ge1$,
$\phi$ is a locally integrable function in ${\mathbb R}^N$, 
$\{a_\alpha\}\subset{\mathbb R}$, 
and $F$ is a continuous function in ${\mathbb R}^N\times[0,\infty)\times{\mathbb R}\times\cdots\times{\mathbb R}^{N^m}$. 
Here $G=G(x,y,t)$ is an integral kernel, which is a generalization of the fundamental solutions 
to the heat equation, fractional heat equations, and higher-order heat equations. 
Throughout this paper we assume the following condition~(G) on the integral kernel~$G$: 
\begin{itemize}
  \item[(G)] 
  \begin{itemize}
  	\item[(a)] 
	$G\in C^{\ell+m}({\mathbb R}^{2N}\times(0,T_*))$ for some $T_*\in(0,\infty]$;
	\item[{\rm (b)}]
	There exist $C_G>0$, $d>\ell+m$, and $L>0$ such that
	$$
	|\nabla^j G(x,y,t)|\le C_G\,t^{-\frac{N+j}{d}}\left(1+t^{-\frac{1}{d}}|x-y|\right)^{-N-L-j}
	$$
	for $(x,y,t)\in{\mathbb R}^{2N}\times(0,T_*)$ and $j\in\{0,\dots,\ell+m\}$;
	\item[(c)] 
	$\displaystyle{G(x,z,t)=\int_{{\mathbb R}^N}G(x,y,t-s)G(y,z,s)dy}\quad$ 
	for $x$, $z\in{\mathbb R}^N$ and $0<s<t<T_*$.
  \end{itemize}
\end{itemize}
The purpose of this paper is to obtain sufficient conditions for
the existence of solutions to integral equation~(I) 
under condition~(G) and a suitable structure condition on~$F$. 
Our sufficient conditions enable us to study the existence of solutions to the Cauchy problem for 
nonlinear parabolic equations of the form
\begin{equation}
\label{eq:1.1}
\partial_t u+{\mathcal L}u=\displaystyle{\sum_{|\alpha|=\ell}}
a_\alpha\partial_x^\alpha F(x,t,u,\dots,\nabla_x^m u).
\end{equation}
Here $-{\mathcal L}$ is a generalization of 
elliptic operators with variable coefficients, fractional elliptic operators, and higher-order elliptic operators. 
\vspace{3pt}

Let us consider the Cauchy problem for the semilinear parabolic equation
\begin{equation}
\tag{S}
\left\{
\begin{array}{ll}
\partial_t u+(-\Delta)^{\frac{d}{2}} u=|u|^p, & \quad x\in{\mathbb R}^N,\,\,t>0,\vspace{3pt}\\
u(x,0)=\phi(x)\ge 0, & \quad x\in{\mathbb R}^N,
\end{array}
\right.
\end{equation}
where $d>0$ and $p>1$. 
The solvability of problem~(S) has been studied in many papers.  
Here we just refer to the monograph~\cite{QS} 
and papers~\cites{BP, C, FL, GP, GG, HI02, HI01, HIT, IKO, KY, P, Q, RS, W, Y},  
which are closely related to this paper. 
Among others, 
for the case of $0<d\le 2$, 
the first author of this paper and Hisa \cite{HI01} developed the arguments in \cites{LN, RS, S} 
and obtained necessary conditions and sufficient conditions for the existence of solutions to problem~(S). 
As corollaries of their main results, they proved the following properties for $0<d\le 2$: 
\begin{itemize}
  \item[(a)] Let $1<p<1+d/N$. Then problem~(S) possesses a local-in-time nonnegative solution 
  if and only if $\sup_{x\in{\mathbb R}^N}\|\phi\|_{L^1(B(x,1))}<\infty$;
  \item[(b)]
  There exists $\gamma>0$ such that, if 
  $$
  \phi(x)\ge
  \left\{
  \begin{array}{ll}
  \gamma|x|^{-\frac{d}{p-1}} & \mbox{if}\quad \displaystyle{p>1+\frac{d}{N}},\vspace{3pt}\\
  \gamma|x|^{-N}\displaystyle{\biggr|\log\biggr(e+\frac{1}{|x|}\biggr)\biggr|^{-\frac{N}{d}-1}} 
   & \mbox{if}\quad \displaystyle{p=1+\frac{d}{N}},\vspace{3pt}\\
  \end{array}
  \right.
  \quad x\in B(0,1),
  $$
  then problem~(S) possesses no local-in-time nonnegative solutions;
  \item[(c)] 
  There exists $\gamma'>0$ such that, if 
  $$
  0\le\phi(x)\le
  \left\{
  \begin{array}{ll}
  \gamma'|x|^{-\frac{d}{p-1}}+C & \mbox{if}\quad \displaystyle{p>1+\frac{d}{N}},\vspace{3pt}\\
  \gamma'|x|^{-N}\displaystyle{\biggr|\log\biggr(e+\frac{1}{|x|}\biggr)\biggr|^{-\frac{N}{d}-1}}+C & \mbox{if}\quad \displaystyle{p=1+\frac{d}{N}},\vspace{3pt}\\
  \end{array}
  \right.
  \quad x\in{\mathbb R}^N,
  $$
  for some $C>0$, 
  then problem~(S) possesses a local-in-time nonnegative solution. 
\end{itemize}
In the proof of assertion~(c), 
it is crucial to construct supersolutions to problem~(S) by the semigroup property of 
the corresponding semigroup.

Subsequently, in \cite{IKO}, the authors of this paper obtained
necessary conditions and sufficient conditions for the existence of solutions to
problem~(S) in the case of $d=2,4,\cdots$,
and proved that assertions~(a), (b), and (c) hold.  
One of the main difficulties in the study of sufficient conditions in the case of $d=2,4,\cdots$ 
comes from the sign-change 
of the fundamental solution~$G_d$ to the parabolic equation
$$
\partial_t u+(-\Delta)^{\frac{d}{2}} u=0,\quad x\in{\mathbb R}^N,\,\,\,t>0.
$$
In order to overcome the difficulty, 
they introduced a majorant kernel $K=K(x,t)$ satisfying 
$$
|G_d(x,t)|\le K(x,t),\quad
\int_{{\mathbb R}^N}K(x-y,t-s)K(y,s)\,dy\le CK(x,t),
$$
for $x\in{\mathbb R}^N$ and $0<s<t$. 
Here $C$ is a positive constant independent of $x\in{\mathbb R}^N$ and $0<s<t$. 
Thanks to the majorant kernel $K$, 
they developed the arguments in \cite{HI01} to obtain sufficient conditions for the existence of solutions to
problem~(S) in the case of $d=2,4,\cdots$. 
\vspace{3pt}

In this paper, under condition~(G), 
we develop the arguments in \cite{IKO} and obtain sufficient conditions 
for  the existence of solutions to integral equation~(I). 
Furthermore, we apply our main results 
to the Cauchy problem for some concrete nonlinear parabolic equations of form~\eqref{eq:1.1}
and obtain rather sharp sufficient conditions for the existence of solutions to the Cauchy problem. 
(See Section~7.) 
\vspace{3pt}

We introduce some notations and formulate the solution to integral equation~(I). 
For any $x\in{\mathbb R}^N$ and $\sigma>0$, let $B(x,\sigma):=\{y\in{\mathbb R}^N\,:\,|y-x|<\sigma\}$ 
and $|B(x,\sigma)|$ denotes the volume of the ball {\bf$B(x,\sigma)$}. 
For any multi-index $\alpha=(\alpha_1,\cdots,\alpha_N)\in(\mathbb N\cup\{0\})^N$, 
we write 
$$
|\alpha|:=\sum_{i=1}^N\alpha_i\quad\mbox{and}\quad
\partial_x^\alpha:=\frac{\partial^{|\alpha|}}{\partial x_1^{\alpha_1}\cdots\partial x_N^{\alpha_N}}.
$$
For any $1\le r\le\infty$, we set $L^r:=L^r({\mathbb R}^N)$ and 
we denote by $L^r_{{\rm uloc}}$ the uniformly local $L^r$ space, that is, 
$f\in L^r_{{\rm uloc}}$ if and only if 
$$
\sup_{x\in{\mathbb R}^N}\|f\|_{L^r(B(x,1))}<\infty. 
$$
Furthermore, for any $1\le q\le r<\infty$, 
we say that a function $f\in L^q_{{\rm uloc}}$ belongs to 
the Morrey space ${\mathcal M}_{r,q}$ if 
$$
\|f\|_{{\mathcal M}_{r,q}}:=\sup_{x\in{\mathbb R}^N}\sup_{\sigma>0}\,
\sigma^{\frac{N}{r}}\left(\,\dashint_{B(x,\sigma)}|\phi(y)|^q\,dy\right)^{\frac{1}{q}}<\infty,
$$
where
$$
\dashint_{B(x,\sigma)} \, f(y)\,dy:=\frac{1}{|B(x,\sigma)|}\int_{B(x,\sigma)} f(y)\,dy.
$$
Set 
$$
D_m:=m+1\quad\mbox{if}\quad N=1,\qquad D_m:=(N^{m+1}-1)/(N-1)\quad\mbox{if}\quad N\ge 2. 
$$
\begin{definition}
\label{Definition:1.1}
Assume condition~{\rm (G)}. Let $F\in C({\mathbb R}^N\times[0,\infty)\times{\mathbb R}^{D_m})$ 
and $0<T\le\infty$. 
We say that $u$ is a solution to integral equation~{\rm (I)} in ${\mathbb R}^N\times [0,T)$ if 
$u\in BC^{m;0}({\mathbb R}^N\times(0,T))$, that is,  
$$
\nabla^j_x u\in BC({\mathbb R}^N\times(\tau,T)),\quad j\in\{0,\dots,m\},
$$
for $\tau\in(0,T)$ and $u$ satisfies~{\rm (I)} for $(x,t)\in{\mathbb R}^N\times(0,T)$.
\end{definition}

We are ready to state our main results. 
Theorem~\ref{Theorem:1.1} is a modification of \cite{IKO}*{Theorem~4.1}
and it is crucial in our study. 
\begin{theorem}
\label{Theorem:1.1}
Let $\ell$, $m\in\{0,1,\dots\}$ 
and let $G$ be the integral kernel satisfying condition~{\rm (G)} for some $L>0$, $d>0$,  
and $T_*\in(0,\infty]$. 
Let $0<\theta<2$ be such that $\theta\le \min\{d,L\}$ and 
set $P_\theta=P_\theta(x,t)$ be the fundamental solution to 
the fractional heat equation
\begin{equation}
\label{eq:1.2}
\partial_tu+(-\Delta)^{\frac{\theta}{2}}u=0,\quad x\in{\mathbb R}^N,\,\,t>0.
\end{equation}
Set
\begin{equation}
\label{eq:1.3}
K_\theta(x,t):=P_\theta\left(x,t^{\frac{\theta}{d}}\right),\quad x\in{\mathbb R}^N,\,\,t>0.
\end{equation}
\begin{itemize}
  \item[{\rm (a)}] 
  For any $j\in\{0,\dots,\ell+m\}$, there exists $c_j>0$ such that 
  \begin{equation}
  \label{eq:1.4}
  |\nabla^j_x G(x,y,t)|\le c_j t^{-\frac{j}{d}}K_\theta(x-y,t),\quad x,y\in{\mathbb R}^N,\,\,t\in(0,T_*).
  \end{equation}
  \item[{\rm (b)}]  
  There exists $C_*>0$ such that 
  $$
  \int_{{\mathbb R}^N}K_\theta(x-y,t-s)K_\theta(y,s)\,dy\le C_*K_\theta(x,t),
  \quad x\in{\mathbb R}^N,\,\,0<s<t.
  $$
\end{itemize}
\end{theorem}
On the basis of Theorem~\ref{Theorem:1.1}, 
we study the existence of solutions to integral equation~(I) 
under the following structure condition~($\mbox{F}_n$) for some $n\in\{0,\dots,m\}$:
\begin{itemize}
  \item[($\mbox{F}_n$)]  
  \begin{itemize}
  	\item[(a)] 
	Let $F$ is a continuous function in ${\mathbb R}^N\times[0,\infty)\times{\mathbb R}^{D_m}$; 
	\item[(b)] 
	There exist $J\subset\{n,\dots,m\}$, ${\bf p}:=\{p_j\}_{j\in J}\subset (0,\infty)$, and $A>-1$ 
	such that $|{\bf p}|:=\displaystyle{\sum_{j\in J}p_j>1}$ and 
	$$
	|F(x,t,z_0,z_1,\dots,z_m)|\le t^A\prod_{j\in J}|z_j|^{p_j}
	$$
	for $(x,t)\in{\mathbb R}^N\times[0,\infty)$ and $z_j\in{\mathbb R}^{N^j}$, where $j\in\{0,\dots,m\}$;
	\item[(c)] Let $\langle {\bf p}\rangle_n:=n+\displaystyle{\sum_{j\in J}} (j-n)p_j$ satisfy
	$$
	d(1+A)\ge\langle {\bf p}\rangle_n+\ell\quad\mbox{if}\quad n+\ell>0,
	\quad
	d(1+A)>\langle {\bf p}\rangle_0\quad\mbox{if}\quad n+\ell=0.
	$$
  \end{itemize}
\end{itemize}
Set 
\begin{equation}
\label{eq:1.5}
r_n:=
\left\{
\begin{array}{ll}
\displaystyle{\frac{N(|{\bf p}|-1)}{d(1+A)-\langle{\bf p}\rangle_n-\ell}} 
& \mbox{if}\quad d(1+A)>\langle {\bf p}\rangle_n+\ell,\vspace{5pt}\\
\infty & \mbox{if}\quad d(1+A)=\langle {\bf p}\rangle_n+\ell.
\end{array}
\right.
\end{equation}
Under structure condition~($\mbox{F}_n$), 
we consider the following four cases:
\begin{equation*}
\begin{array}{ll}
{\rm (A)}\quad 0<n+\ell<d(1+A); \qquad& 
{\rm (B)}\quad\mbox{$n=\ell=0$ and $r_0<1$};\vspace{7pt}\\
{\rm (C)}\quad\mbox{$n=\ell=0$ and $r_0>1$};\,\,\, & 
{\rm (D)}\quad\mbox{$n=\ell=0$ and $r_0=1$},
\end{array}
\end{equation*}
and state our sufficient conditions for the existence of solutions to integral equation~(I). 
Our sufficient conditions are represented in the spirit of Morrey spaces and their generalizations. 

\begin{theorem}
\label{Theorem:1.2}
Assume conditions~{\rm (G)} and {\rm ($\mbox{F}_n$)} for some $n\in\{0,\dots,m\}$.  
Consider case~{\rm (A)}. 
Let $\phi\in L^1_{{\rm uloc}}$. 
Then there exists $\gamma>0$ such that, if
\begin{equation}
\label{eq:1.6}
\sup_{x\in{\mathbb R}^N}\,\sup_{0<\sigma<T^{\frac{1}{d}}}\,
\sigma^{\frac{N}{r_n}}\,\dashint_{B(x,\sigma)}|\nabla^n \phi(y)|\,dy\le\gamma
\end{equation}
for some $T\in(0,T_*]$,
integral equation~{\rm (I)} possesses a solution $u$ in ${\mathbb R}^N\times[0,T)$ such that 
\begin{equation}
\label{eq:1.7}
|\nabla^j u(x,t)|
\le
\left\{
\begin{array}{ll}
 CT^{-\frac{N}{d}\left(\frac{1}{r_n}-1\right)}t^{-\frac{N}{d}-\frac{j-n}{d}} & \mbox{if}\quad r_n<1,\vspace{5pt}\\
 Ct^{-\frac{N}{dr_n}-\frac{j-n}{d}} & \mbox{if}\quad r_n\ge 1,
\end{array}
\right.
\end{equation}
for $(x,t)\in{\mathbb R}^N\times(0,T)$ and $j\in\{n,\dots,m\}$. 
Here $C$ is a positive constant independent of~$T$.
\end{theorem}
As a corollary of Theorem~\ref{Theorem:1.2}, we have:
\begin{corollary}
\label{Corollary:1.1}
Assume conditions~{\rm (G)} with $T_*=\infty$ and {\rm ($\mbox{F}_n$)} for some $n\in\{0,\dots,m\}$. 
Consider case~{\rm (A)} and assume $r_n\ge 1$. 
Then there exists $\gamma>0$ such that, if 
$$
\|\nabla^n\phi\|_{{\mathcal M}_{r_n,1}}\le\gamma,
$$ 
integral equation~{\rm (I)} possesses a global-in-time solution~$u$ satisfying 
$$
|\nabla^j u(x,t)|
\le Ct^{-\frac{N}{dr_n}-\frac{j-n}{d}}
$$
for $(x,t)\in{\mathbb R}^N\times(0,\infty)$ and $j\in\{n,\dots,m\}$, where $C$ is a positive constant. 
\end{corollary}
In the following two theorems we treat cases (B) and (C).
\begin{theorem}
\label{Theorem:1.3}
Assume conditions~{\rm (G)} and {\rm ($\mbox{F}_0$)}. 
Consider case {\rm (B)}. 
Then there exists $\gamma>0$ such that, if
\begin{equation}
\label{eq:1.8}
\sup_{x\in{\mathbb R}^N}\dashint_{B(x,T^{\frac{1}{d}})}|\phi(y)|\,dy\le\gamma 
T^{-\frac{N}{dr_0}}
\end{equation}
for some $T\in(0,\infty)$ with $T\le T_*$, 
integral equation~{\rm (I)} possesses a solution~$u$ in ${\mathbb R}^N\times[0,T)$ such that 
$$
|\nabla u(x,t)|\le Ct^{-\frac{N}{d}-\frac{j}{d}}
$$
for $(x,t)\in{\mathbb R}^N\times(0,T)$ and $j\in\{0,\dots,m\}$, where $C$ is a positive constant independent of~$T$.
\end{theorem}
\begin{theorem}
\label{Theorem:1.4}
Assume conditions~{\rm (G)} and {\rm ($\mbox{F}_0$)}. 
Consider case {\rm (C)}.   
Then, for any $q>1$,  there exists $\gamma>0$ such that, if
\begin{equation}
\label{eq:1.9}
\sup_{x\in{\mathbb R}^N}\,\sup_{0<\sigma<T^{\frac{1}{d}}}\,
\sigma^{\frac{N}{r_0}}\,\biggr(\,\dashint_{B(x,\sigma)}|\phi(y)|^q\,dy\biggr)^{\frac{1}{q}}\le\gamma
\end{equation}
for some $T\in(0,T_*]$, integral equation~{\rm (I)} possesses a solution~$u$ in ${\mathbb R}^N\times[0,T)$ such that 
\begin{equation}
\label{eq:1.10}
|\nabla u(x,t)|\le Ct^{-\frac{N}{dr_0}-\frac{j}{d}}
\end{equation}
for $(x,t)\in{\mathbb R}^N\times(0,T)$ and $j\in\{0,\dots,m\}$, 
where $C$ is a positive constant independent of~$T$.
\end{theorem}
As a corollary of Theorem~\ref{Theorem:1.4}, we have:
\begin{corollary}
\label{Corollary:1.2}
Assume conditions~{\rm (G)} with $T_*=\infty$ and {\rm ($\mbox{F}_0$)}. 
Consider case {\rm (C)}. 
Then, for any $q>1$, 
there exists $\gamma>0$ such that, if 
$$
\|\phi\|_{{\mathcal M}_{r_0,q}}\le\gamma,
$$ 
integral equation~{\rm (I)} possesses a global-in-time solution~$u$ 
satisfying \eqref{eq:1.10} in ${\mathbb R}^N\times(0,\infty)$.
\end{corollary}
\begin{remark}
\label{Remark:1.1}
Corollary~{\rm\ref{Corollary:1.2}} is a generalization of \cite{IKK01}*{Theorem~1.1}. 
Indeed, 
assume condition~{\rm ($\mbox{F}_0$)} and consider case~{\rm (C)}. 
Under stronger assumptions than condition~{\rm (G)} with $T_*=\infty$,   
the existence of global-in-time solutions to integral equation~{\rm (I)}
was proved in~\cite{IKK01}*{Theorem~1.1} for the case 
when $\phi\in W^{1,m}({\mathbb R}^N)$ and $\|\phi\|_{L^{r_0,\infty}}$ is small enough. 
Here $L^{r_0,\infty}$ is the weak $L^{r_0}$ space in ${\mathbb R}^N$. 
On the other hand, 
by Corollary~{\rm\ref{Corollary:1.2}} we easily obtain 
the existence of global-in-time solutions to integral equation~{\rm (I)} 
provided that $\|\phi\|_{L^{r_0,\infty}}$ is small enough, since 
$L^{r_0,\infty}\subset{{\mathcal M}_{r_0,q}}$ for $1\le q<r_0$ 
{\rm ({\it see} \cite{KY}*{Lemma~1.7})}. 
\end{remark}
We state our result in case (D).
\begin{theorem}
\label{Theorem:1.5}
Assume conditions~{\rm (G)} and {\rm ($\mbox{F}_0$)}. 
Consider case {\rm (D)}. Let $\beta>0$. Set 
\begin{equation}
\label{eq:1.11}
\Phi(s):=s[\log(e+s)]^\beta,
\quad
\rho(s):=s^{-N}\left[\log\left(e+s^{-1}\right)\right]^{-\frac{N}{d(1+A)-\langle{\bf p}\rangle_0}}.
\end{equation}
Then there exists $\gamma>0$ such that, if 
\begin{equation}
\label{eq:1.12}
\sup_{x\in{\mathbb R}^N}\Phi^{-1}\biggr(\,\dashint_{B(x,\sigma)}\Phi(T^{\frac{N}{d}}|\phi(y)|)\,dy\biggr)\le\gamma \rho(\sigma T^{-\frac{1}{d}}),
\quad 0<\sigma<T^{\frac{1}{d}},
\end{equation}
for some $T\in(0,\infty)$ with $T\le T_*$, integral equation~{\rm (I)} possesses a solution~$u$ in ${\mathbb R}^N\times[0,T)$ such that 
\begin{equation}
\label{eq:1.13}
|\nabla^j u(x,t)|
\le Ct^{-\frac{N}{d}-\frac{j}{d}}\biggr|\log\left(\frac{t}{2T}\right)\biggr|^{-\frac{N}{d}},\quad j\in\{0,\dots,m\},
\end{equation}
for $(x,t)\in{\mathbb R}^N\times(0,T)$, where $C$ is a positive constant independent of $T$.
\end{theorem}
As a corollary of Theorem~\ref{Theorem:1.5}, we have:
\begin{corollary}
\label{Corollary:1.3}
Assume conditions~{\rm (G)} and {\rm ($\mbox{F}_0$)}. 
Consider case {\rm (D)}. 
Then there exists $\gamma>0$ such that, if 
$$
\sup_{x\in{\mathbb R}^N}\int_{B(x,\sigma)} |\phi(y)|\left[\log\left(e+T^{\frac{N}{d}}|\phi(y)|\right)\right]
^{\frac{N}{d(1+A)-\langle{\bf p}\rangle_0}}\,dy
\le\gamma ,
\quad 0<\sigma<T^{\frac{1}{d}},
$$
for some $T\in(0,\infty)$ with $T\le T_*$, integral equation~{\rm (I)} possesses a solution~$u$ in ${\mathbb R}^N\times[0,T)$ 
satisfying \eqref{eq:1.13} in ${\mathbb R}^N\times(0,T)$.
\end{corollary}

Here we mention the strategy for the proofs of our sufficient conditions. 
As the first step we construct approximate solutions~$\{u_\epsilon\}_{\epsilon>0}$ to integral equation~(I). 
Next, thanks to the integral kernel $K_\theta$ given in Theorem~\ref{Theorem:1.1}, 
we develop the arguments in \cite{IKO} to find supersolutions to integral equation~(I). 
This enables us to obtain uniform estimates of approximate solutions~$\{u_\epsilon\}_{\epsilon>0}$. 
Then, applying the parabolic regularity theorems and the Arzel\'a--Ascoli theorem,  
we find a solution to integral equation~(I), and the proofs of 
our sufficient conditions are complete. 

The rest of this paper is organized as follows: 
In Section~2 we recall some properties of fundamental solutions to fractional heat equations 
and prove Theorem~\ref{Theorem:1.1}. 
In Section~3 we construct approximate solutions to integral equation~(I) 
and obtain some decay estimates of the approximate solutions. 
In Sections~4, 5, and 6 we find supersolutions to integral equation~(I) 
and prove our main theorems and their corollaries.  
In Section~7 we apply our main theorems to the Cauchy problem for some concrete nonlinear parabolic equations 
and we show the sharpness of our sufficient conditions. 
\section{Proof of Theorem~\ref{Theorem:1.1}}
%
In this section we prove Theorem~\ref{Theorem:1.1}. 
In what follows, 
by $C$ we denotes generic positive constants 
and they may have different values also within the same line.  
We recall some properties of the fundamental solution~$P_\theta$ to the fractional 
heat equation~\eqref{eq:1.2}, 
where $0<\theta<2$. 
The fundamental solution~$P_\theta$ is represented by 
$$
P_\theta(x,t)=(2\pi)^{-\frac{N}{2}}\int_{{\mathbb R}^N}e^{ix\cdot \xi}e^{-t|\xi|^\theta}\,d\xi.
$$
Then $P_\theta=P_\theta(x,t)$ is a positive, smooth, and radially symmetric function 
in ${\mathbb R}^N\times(0,\infty)$ 
and satisfies the following properties (see \cites{BJ,BK}): 
\begin{align}
\label{eq:2.1}
& P_\theta(x,t)=t^{-\frac{N}{\theta}}P_\theta(t^{-\frac{1}{\theta}}x,1),\\
\label{eq:2.2}
 &|(\nabla^j P_\theta)(x,t)|\le C_j t^{-\frac{N+j}{\theta}}
\big(1+t^{-\frac{1}{\theta}}|x|\big)^{-N-\theta-j},\\
 \label{eq:2.3}
& P_\theta(x,t)\ge 
C^{-1} t^{-\frac{N}{\theta}}
\big(1+t^{-\frac{1}{\theta}}|x|\big)^{-N-\theta},\\
\label{eq:2.4}
&
\int_{\mathbb R^N}P_\theta(x,t)\,dx=1,
\end{align}
for $x\in{\mathbb R}^N$, $t>0$, and $j=0,1,2,\dots$. 
Here $C_j$ is a positive constant depending on $j$.
Furthermore,
\begin{equation}
\label{eq:2.5}
P_\theta(x,t)=\int_{{\mathbb R}^N}P_\theta(x-y,t-s)\,P_\theta(y,s)\,dy,
\quad x\in{\mathbb R}^N,\,\,0<s<t.
\end{equation}
For any $\phi\in L^1_{\rm uloc}$, 
we set 
\begin{equation}
\label{eq:2.6}
[S_\theta(t)\phi](x):=\int_{{\mathbb R}^N}P_\theta(x-y,t)\phi(y)\,dy.
\end{equation}
Then, for any $j=0,1,2,\dots$, 
by the Young inequality and \eqref{eq:2.2}
we find $C_j'>0$ such that 
$$
\|\nabla^j S_\theta(t)\phi\|_{L^q}\le C_j'
t^{-\frac{N}{\theta}(\frac{1}{p}-\frac{1}{q})-\frac{j}{\theta}}\|\phi\|_{L^p},\quad t>0,
$$
for $\phi\in L^q$ and $1\le p\le q\le\infty$.  
(See e.g. \cite{IKK01}*{Section~2}.)  
Furthermore, we recall the following lemma on the decay of $\|S_\theta(t)\phi\|_{L^\infty}$ (see \cite{HI01}*{Lemma~2.1}). 
\begin{lemma}
\label{Lemma:2.1}
Let $0<\theta< 2$. 
Then there exists $C=C(N,\theta)>0$ such that
$$
 \|S_\theta(t)\mu\|_{L^\infty}
 \le Ct^{-\frac{N}{\theta}}\sup_{x\in{\mathbb R}^N}|\mu(B(x,t^{\frac{1}{\theta}}))|,
 \quad t>0,
$$
for Radon measures $\mu$ in ${\mathbb R}^N$.
\end{lemma}

We prove Theorem~\ref{Theorem:1.1}. 
\vspace{3pt}
\newline
{\bf Proof of Theorem~\ref{Theorem:1.1}.} 
For any $j\in\{0,\dots,\ell+m\}$,
it follows from condition~(G), \eqref{eq:1.3}, and \eqref{eq:2.1} that
\begin{equation*}
\begin{split}
|\nabla_x^j G(x,y,t)| & 
\le Ct^{-\frac{N+j}{d}}\left(1+t^{-\frac{1}{d}}|x-y|\right)^{-N-L-j}\\
&
\le Ct^{-\frac{N+j}{d}}\left(1+t^{-\frac{1}{d}}|x-y|\right)^{-N-\theta}\\
 & \le Ct^{-\frac{N+j}{d}}P_\theta\left(t^{-\frac{1}{d}}(x-y),1\right)
 =Ct^{-\frac{j}{d}}K_\theta(x-y,t)
 \end{split}
\end{equation*}
for $x, y\in{\mathbb R}^N$ and $0<t<T_*$. 
Here we used the assumption that $\theta\le L$. 
This implies assertion~(a). 
On the other hand, 
since $\theta\le d$, we have
$$
t^{\frac{\theta}{d}}=(t-s+s)^{\frac{\theta}{d}}
\le(t-s)^{\frac{\theta}{d}}+s^{\frac{\theta}{d}}\le 2t^{\frac{\theta}{d}}
$$
for $0<s<t$.
Then, by \eqref{eq:1.3} and \eqref{eq:2.5} we have
\begin{equation*}
\begin{split}
 & \int_{{\mathbb R}^N}K_\theta(x-y,t-s)K_\theta(y,s)\,dy\\
 & =\int_{{\mathbb R}^N}P_\theta\left(x-y,(t-s)^{\frac{\theta}{d}}\right)
 P_\theta\left(y,s^{\frac{\theta}{d}}\right)\,dy
 \\
 &
=P_\theta\left(x,(t-s)^{\frac{\theta}{d}}+s^{\frac{\theta}{d}}\right)\\
 & =P_\theta\left(x,\kappa_{t,s}t^{\frac{\theta}{d}}\right),
 \quad\mbox{where}\quad
 \kappa_{t,s}:=\frac{{(t-s)^{\frac{\theta}{d}}+s^{\frac{\theta}{d}}}}{t^{\frac{\theta}{d}}}\in[1,2],
\end{split}
\end{equation*}
for $x\in{\mathbb R}^N$ and $0<s<t$.
Furthermore, it follows from \eqref{eq:2.2} and \eqref{eq:2.3} that
\begin{equation*}
\begin{split}
P_\theta(x,\kappa_{t,s} t) 
& 
\le C(\kappa_{t,s}t)^{-\frac{N}{\theta}}\left(1+(\kappa_{t,s} t)^{-\frac{1}{\theta}}|x|\right)^{-N-\theta}\\
 & \le Ct^{-\frac{N}{\theta}}\left(1+t^{-\frac{1}{\theta}}|x|\right)^{-N-\theta}
\le CP_\theta(x,t)
\end{split}
\end{equation*}
for $x\in{\mathbb R}^N$ and $0<s<t$. 
These imply assertion~(b). 
Thus Theorem~\ref{Theorem:1.1} follows.
$\Box$\vspace{5pt}

Similarly to \eqref{eq:2.6}, 
we set 
$$
[S(t)\phi](x):=\int_{{\mathbb R}^N}G(x,y,t)\phi(y)\,dy,
\quad
[S_{K_\theta}(t)\phi](x):=\int_{{\mathbb R}^N}K_\theta(x-y,t)\phi(y)\,dy,
$$
for $\phi\in L^1_{\rm uloc}$.
Then we observe from Theorem~\ref{Theorem:1.1} that 
\begin{align}
\label{eq:2.7}
 & |[\nabla^j S(t)\phi](x)|\le c_jt^{-\frac{j}{d}}[S_{K_\theta}(t)|\phi|](x),
\quad t\in(0,T_*),
\quad j\in\{0,\dots, \ell+m\},\\
\label{eq:2.8}
 & [S_{K_\theta}(t-s)[S_{K_\theta}(s)\phi]](x)\le C_*[S_{K_\theta}(t)\phi](x),
 \quad 0<s<t,
\end{align}
for $x\in{\mathbb R}^N$. 
Furthermore, it follows from Lemma~\ref{Lemma:2.1} with \eqref{eq:1.3} that
\begin{equation}
\label{eq:2.9}
 \|S_{K_\theta}(t)\phi\|_{L^\infty}
 \le C\sup_{x\in{\mathbb R}^N}\dashint_{B(x,t^{\frac{1}{d}})}|\phi(y)|\,dy,
 \quad t>0.
\end{equation}
These properties are crucial in the proof of our sufficient conditions 
for the existence of solutions to integral equation~(I). 
\section{Approximate solutions}
Let $\ell$, $m\in\{0,1,\dots\}$. Assume condition~($\mbox{F}_n$) for some $n\in\{0,\dots,m\}$. 
We construct approximate solutions to integral equation~(I).
For $\epsilon>0$, let
$$
F_\epsilon(x,t,z)
:=
\left\{
\begin{array}{ll}
-\epsilon^{-1}\quad & \mbox{if}\quad F(x,|t|,z)<-\epsilon^{-1},\vspace{3pt}\\
F(x,|t|,z)\quad & \mbox{if}\quad -\epsilon^{-1}\le F(x,|t|,z)\le\epsilon^{-1},\vspace{3pt}\\
\epsilon^{-1}\quad & \mbox{if}\quad F(x,|t|,z)>\epsilon^{-1},\vspace{3pt}
\end{array}
\right.
$$
for $(x,t)\in{\mathbb R}^{N+1}$ and $z=(z_1,\dots,z_m)\in {\mathbb R}^{D_m}$.
Let $\rho\in C_0^\infty({\mathbb R}^{N+1+D_m})$ be such that 
$$
\rho\ge 0\quad\mbox{in}\quad {\mathbb R}^{N+1+D_m},
\quad
\rho=0\quad\mbox{if}\quad |(x,t,z)|\ge 1,
\quad
\int_{{\mathbb R}^{N+1+D_m}}\rho(x,t,z)\,dx\,dt\,dz=1.
$$
Set 
$$
\tilde{F}_\epsilon(x,t,z):=\epsilon^{-N-1-D_m}
\int_{{\mathbb R}^{N+1+D_m}}
\rho\left(\epsilon(x-y),\epsilon(t-s),\epsilon(z-\xi)\right)F_\epsilon(y,s,\xi)\,dy\,ds\,d\xi
$$
for $(x,t,z)\in{\mathbb R}^{N+1+D_m}$.
Then we easily see that
\begin{equation}
\label{eq:3.1}
\tilde{F}_\epsilon\in BC^m({\mathbb R}^{N+1+D_m}),
\quad  \|\tilde{F}_\epsilon\|_{L^\infty({\mathbb R}^{N+1+D_m})}\le\epsilon^{-1},
\quad  \|\nabla_z\tilde{F}_\epsilon\|_{L^\infty({\mathbb R}^{N+1+D_m})}\le C\epsilon^{-2}.
\end{equation}
Furthermore, it follows from condition~($\mbox{F}_n$) that 
\begin{align}
\notag
 & \lim_{\epsilon\to +0}\tilde{F}_\epsilon(x,t,z)=F(x,|t|,z),\\
\label{eq:3.2}
 & |\tilde{F}_\epsilon(x,t,z)|
 \le |t|^A\prod_{j\in J}\left(|z_j|+\epsilon\right)^{p_j}
\le C|t|^A\prod_{j\in J}\max\left\{|z_j|,\epsilon\right\}^{p_j},
\end{align}
for $(x,t,z)\in{\mathbb R}^{N+1+D_m}$. 
In this section we prove the following lemma.
\begin{lemma}
\label{Lemma:3.1}
Assume conditions~{\rm (G)} and {\rm($\mbox{F}_n$)} for some $n\in\{0,\dots,m\}$.
Let $\epsilon>0$ and $\tilde{F}_\epsilon$ be as in the above. 
Assume that $\phi$ is a measurable function in ${\mathbb R}^N$ such that 
$$
\sup_{x\in{\mathbb R}^N}\int_{B(x,1)}|\phi(y)|\,dy<\infty.
$$
\begin{itemize}
  \item[{\rm (a)}] 
  There exists $u^\epsilon\in C^{m;0}({\mathbb R}^N\times(0,T_*))$ such that 
   \begin{equation}
  \label{eq:3.3}
   u^\epsilon(x,t)=[S(t)\phi](x) 
  +\sum_{|\alpha|=\ell}a_\alpha\int_0^t\partial_x^\alpha
  \left[S(t-s)\tilde{F}_\epsilon(s,u^\epsilon(s),\dots,\nabla^m u^\epsilon(s))\right](x)\, ds
  \end{equation}
  for $(x,t)\in{\mathbb R}^N\times(0,T_*)$. 
  Here $T_*$ is as in condition~{\rm (G)}. 
  \item[{\rm (b)}] 
  There exists $c_*>0$ with the following property: 
  If there exist $T\in(0,T_*]$ and a continuous function $U^\epsilon$ in ${\mathbb R}^N\times(0,T)$ such that
  \begin{align}
  \label{eq:3.4}
  & t^{-\frac{j-n}{d}}U^\epsilon(x,t)\ge\epsilon,\\
  \label{eq:3.5}
  & c_*[S_{K_\theta}(t)|\nabla^n\phi|](x)
  \le \frac{1}{2}U^\epsilon(x,t),\\
  \label{eq:3.6} 
  & c_*\int_0^t(t-s)^{-\frac{\ell+j}{d}}s^{A-\frac{\langle {\bf p}\rangle_n-n}{d}}
  \left[S_{K_\theta}(t-s)U^\epsilon(s)^{|{\bf p}|}\right](x)\,ds
  \le\frac{1}{2}t^{-\frac{j-n}{d}}U^\epsilon(x,t),
  \end{align}
  for $(x,t)\in{\mathbb R}^N\times(0,T)$ and $j\in\{n,\dots,m\}$, then 
  \begin{equation}
  \label{eq:3.7}
  |\nabla^j u^\epsilon(x,t)|\le t^{-\frac{j-n}{d}}U^\epsilon(x,t)
  \end{equation}
  for $(x,t)\in{\mathbb R}^N\times(0,T)$ and $j\in\{n,\dots,m\}$. 
\end{itemize}
\end{lemma}
{\bf Proof.}
Set $u^\epsilon_0(x,t):=[S(t)\phi](x)$ for $(x,t)\in{\mathbb R}^N\times(0,T_*)$. 
It follows from \eqref{eq:2.7} and \eqref{eq:2.9} that 
\begin{equation}
\label{eq:3.8}
\begin{split}
&\qquad
u^\epsilon_0\in C^{m;0}({\mathbb R}^N\times(0,T_*)),
\\
&
\sup_{(x,t)\in{\mathbb R}^N\times[\tau,T_*)}|\nabla^j u^\epsilon_0(x,t)|<\infty,
\quad
\quad\tau\in(0,T_*),\quad  j\in\{0,\dots,\ell+m\}.
\end{split}
\end{equation}
Since $\ell+m<d$, by \eqref{eq:2.7} 
we can define $\{u^\epsilon_k\}\subset BC^{m;0}({\mathbb R}^N\times(0,T_*))$ inductively
as follows:
\begin{equation}
\label{eq:3.9}
u^\epsilon_{k+1}(x,t):=u^\epsilon_0(x,t)
+\sum_{|\alpha|=\ell}a_\alpha\int_0^t\partial_x^\alpha
  \left[S(t-s)\tilde{F}_\epsilon(s,u^\epsilon_k(s),\dots,\nabla^m u^\epsilon_k(s))\right](x)\, ds
\end{equation}
for $(x,t)\in{\mathbb R}^N\times(0,T_*)$ and $k=0,1,2,\dots$. 
Then it follows from condition~(G)~(c) that
\begin{equation}
\label{eq:3.10}
u^\epsilon_{k+1}(x,t):=u^\epsilon_{k+1}(x,\tau)
+\sum_{|\alpha|=\ell}a_\alpha\int_\tau^t\partial_x^\alpha
  \left[S(t-s)\tilde{F}_\epsilon(s,u^\epsilon_k(s),\dots,\nabla^m u^\epsilon_k(s))\right](x)\, ds
\end{equation}
for $(x,t)\in{\mathbb R}^N\times(\tau,T_*)$, $\tau\in(0,T_*)$, and $k=0,1,2,\dots$. 

Let $\delta\in(0,T_*)$ be small enough. 
Set 
$$
L_k:=\sup_{j\in\{0,\dots,m\}}\sup_{(x,t)\in{\mathbb R}^N\times(0,\delta]}|\nabla^j u^\epsilon_{k+1}(x,t)-\nabla^j u^\epsilon_k(x,t)|.
$$
Thanks to the mean value theorem, 
by \eqref{eq:1.4}, \eqref{eq:3.1}, and \eqref{eq:3.9}
we obtain 
\begin{equation*}
\begin{split}
 & |\nabla^j u_{k+1}^\epsilon(x,t)-\nabla^j u_k^\epsilon(x,t)|\\
 & \le CL_{k-1}\int_0^t\int_{{\mathbb R}^N}|\nabla^{\ell+j}G(x,y,t-s)|
 \|\nabla_z \tilde F_\epsilon\|_{L^\infty({\mathbb R}^{N+1+D_m})}\,dy\,ds\\
 & \le C\epsilon^{-2}L_{k-1}\int_0^t\int_{{\mathbb R}^N}(t-s)^{-\frac{\ell+j}{d}}K_\theta(x-y,t-s)\,dy\,ds\\
 & \le C\epsilon^{-2}L_{k-1}\delta^{1-\frac{\ell+j}{d}}
\end{split}
\end{equation*}
for $(x,t)\in{\mathbb R}^N\times(0,\delta]$, $j\in\{0,\dots,m\}$, and $k=1,2,\dots$. 
Then, taking small enough $\delta>0$ if necessary, we obtain
$$
L_k\le C\epsilon^{-2}L_{k-1}\delta^{1-\frac{\ell+j}{d}}\le\frac{1}{2}L_{k-1},\quad k=1,2,\dots.
$$
This together with \eqref{eq:3.8} implies that 
$\{u^\epsilon_k\}_{k=0}^\infty$ is a Cauchy sequence in $BC^{m;0}({\mathbb R}^N\times(\tau,\delta])$ for $\tau\in(0,\delta)$.
Therefore there exists $u^\epsilon\in C^{m;0}({\mathbb R}^N\times(0,\delta])$ such that 
\begin{equation}
\label{eq:3.11}
u^\epsilon\in BC^{m;0}({\mathbb R}^N\times(\tau,\delta]),
\quad
\lim_{k\to\infty}\sup_{(x,t)\in{\mathbb R}^N\times(\tau,\delta]}|\nabla^j u^\epsilon_k(x,t)-\nabla^j u^\epsilon(x,t)|=0
\end{equation}
for $\tau\in(0,\delta)$. 
Then we apply the Lebesgue dominated convergence theorem to see that 
$u^\epsilon$ satisfies \eqref{eq:3.3} in ${\mathbb R}^N\times(0,\delta]$. 
This implies that assertion~(a) holds for $t\in(0,\delta]$.
Repeating the above arguments with $\phi$ replaced by $u^\epsilon(\cdot,\delta)$,
due to \eqref{eq:3.10}, 
we see that assertion~(a) holds for $t\in(0,2\delta]\cap(0,T_*)$.
Therefore, repeating this argument several times, we deduce that assertion~(a) holds. 

We prove assertion~(b). 
Set
$$
c_*:=\biggr(1+\sum_{|\alpha|=\ell}|a_\alpha|\biggr)\max_j c_j,
$$
where $c_j$ is as in Theorem~\ref{Theorem:1.1}. 
Assume \eqref{eq:3.4}, \eqref{eq:3.5}, and \eqref{eq:3.6}. 
It follows from \eqref{eq:3.5} that 
\begin{equation}
\label{eq:3.12}
\begin{split}
|\nabla^j u_0^\epsilon(x,t)| 
& \le \int_{{\mathbb R}^N}|\nabla^{j-n} G(x,y,t)||\nabla^n\phi(y)|\,dy\\
 & \le c_{j-n}t^{-\frac{j-n}{d}}\int_{{\mathbb R}^N}K_\theta(x-y,t)|\nabla^n\phi(y)|\,dy
 \le\frac{1}{2}t^{-\frac{j-n}{d}}U^\epsilon(x,t)
\end{split}
\end{equation}
for $(x,t)\in{\mathbb R}^N\times(0,T)$ and $j\in\{n,\dots,m\}$. 
Furthermore, if 
$$
|\nabla^j u_k^\epsilon(x,t)|\le t^{-\frac{j-n}{d}}U^\epsilon(x,t),
\quad(x,t)\in{\mathbb R}^N\times(0,T),\,\,j=n,\dots,m,
$$
for some $k\in\{0,1,2,\dots\}$, then, 
by \eqref{eq:1.4}, \eqref{eq:3.2}, \eqref{eq:3.4}, \eqref{eq:3.6}, \eqref{eq:3.9},
and \eqref{eq:3.12} we find $C>0$ such that
\begin{equation*}
\begin{split}
 & |\nabla^j u^\epsilon_{k+1}(x,t)|
 \le\frac{1}{2}t^{-\frac{j-n}{d}}U^\epsilon(x,t)\\
 & \quad
 +\sum_{|\alpha|=\ell}|a_\alpha|\int_0^t\int_{{\mathbb R}^N}|\nabla^{\ell+j} G(x,y,t-s)|\tilde{F}_\epsilon(y,s,u^\epsilon_k(y,s),\dots,\nabla^m u^\epsilon_k(y,s))|\,dy\,ds\\
 & \le\frac{1}{2}t^{-\frac{j-n}{d}}U^\epsilon(x,t)\\
 & \quad
 +c_{\ell+j}\sum_{|\alpha|=\ell}|a_\alpha|
 \int_0^t\int_{{\mathbb R}^N}(t-s)^{-\frac{\ell+j}{d}}K_\theta(x-y,t-s)s^A\prod_{j\in J}
 \max\{s^{-\frac{j-n}{d}}U^\epsilon(y,s),\epsilon\}^{p_j}\,dy\,ds\\
 & \le\frac{1}{2}t^{-\frac{j-n}{d}}U^\epsilon(x,t)
  +c_*\int_0^t\int_{{\mathbb R}^N}(t-s)^{-\frac{\ell+j}{d}}K_\theta(x-y,t-s)s^A\prod_{j\in J}
 \left(s^{-\frac{j-n}{d}}U^\epsilon(y,s)\right)^{p_j}\,dy\,ds\\
  & \le\frac{1}{2}t^{-\frac{j-n}{d}}U^\epsilon(x,t)
 +c_*\int_0^t\int_{{\mathbb R}^N}(t-s)^{-\frac{\ell+j}{d}}
 s^{A-\frac{\langle {\bf p}\rangle_n-n}{d}}K_\theta(x-y,t-s)U^\epsilon(y,s)^{|{\bf p}|}\,dy\,ds\\
 & \le t^{-\frac{j-n}{d}}U^\epsilon(x,t)
\end{split}
\end{equation*}
for $(x,t)\in{\mathbb R}^N\times(0,T)$. 
This together with \eqref{eq:3.12} implies that 
$$
|\nabla^j u^\epsilon_k(x,t)|\le t^{-\frac{j-n}{d}}U^\epsilon(x,t)
$$
for $(x,t)\in{\mathbb R}^N\times(0,T)$, $j\in\{n,\dots,m\}$, and $k=0,1,2,\dots$. 
Then, thanks to \eqref{eq:3.11}, we obtain \eqref{eq:3.7}. 
Thus assertion~(b) holds, 
and Lemma~\ref{Lemma:3.1} follows. 
$\Box$
\section{Proof of Theorem~\ref{Theorem:1.2}}
In this section we prove Theorem~\ref{Theorem:1.2} and Corollary~\ref{Corollary:1.1}.
\vspace{3pt}
\newline
{\bf Proof of Theorem~\ref{Theorem:1.2}.}
Let $T\in(0,T_*]$ and assume \eqref{eq:1.6}. 
Let $c_*$ be as in Lemma~\ref{Lemma:3.1}. 
Let $i\in\{1,2,\dots\}$ and fix it. 
Let $T_i:=\min\{T,i\}$ and $\epsilon_i\in(0,1)$ be small enough. 
Then we find $L_i>0$ such that 
$$
c_*L_i\min_{j\in\{n,\dots,m\}}T_i^{-\frac{j-n}{d}}=\epsilon_i.
$$ 
Set
\begin{equation}
\label{eq:4.1}
U^i(x,t):=2c_*\left[S_{K_\theta}(t)\left(|\nabla^n\phi|+L_i\right)\right](x)
=2c_*\left[S_{K_\theta}(t)|\nabla^n\phi|\right](x)+2c_*L_i.
\end{equation}
Then we see that 
\begin{equation}
\label{eq:4.2}
\begin{split}
 & \inf_{(x,t)\in{\mathbb R}^N\times(0,T_i)}t^{-\frac{j-n}{d}}U^i(x,t)\ge 2c_*L_i T_i^{-\frac{j-n}{d}}\ge\epsilon_i,
 \quad j\in\{n,\dots,m\},\\
 & c_*[S_{K_\theta}(t)|\nabla^n\phi|](x)
 \le \frac{1}{2}U^i(x,t),\quad (x,t)\in{\mathbb R}^N\times(0,T_i).
\end{split}
\end{equation}
Furthermore, by \eqref{eq:1.6} we have
$$
\sup_{x\in{\mathbb R}^N}\int_{B(x,t^{\frac{1}{d}})}|\nabla^n \phi(y)|\,dy
\le\left\{
\begin{array}{ll}
\gamma {T_i}^{\frac{N}{d}\left(1-\frac{1}{r_n}\right)} & \quad\mbox{if}\quad r_n<1,\vspace{3pt}\\
\gamma t^{\frac{N}{d}-\frac{N}{dr_n}} & \quad\mbox{if}\quad r_n\ge 1,
\end{array}
\right.
$$
for $0<t<T_i$. 
Since $\lim_{\epsilon_i\to 0}L_i=0$, 
taking small enough $\epsilon_i>0$ if necessary, 
by \eqref{eq:2.9}
we see that
\begin{equation}
\label{eq:4.3}
U^i(x,t)\le C\gamma_i t^{-\kappa}+2c_*L_i
\le 2\gamma_i Ct^{-\kappa},
\quad (x,t)\in{\mathbb R}^N\times(0,T_i),
\end{equation}
where 
\begin{equation}
\label{eq:4.4}
\gamma_i:=
\left\{
\begin{array}{ll}
\gamma {T_i}^{\frac{N}{d}\left(1-\frac{1}{r_n}\right)} & \mbox{if}\quad r_n<1,\vspace{3pt}\\
\gamma & \mbox{if}\quad r_n\ge 1,
\end{array}
\right.
\qquad
\kappa:=
\left\{
\begin{array}{ll}
\displaystyle{\frac{N}{d}} & \mbox{if}\quad r_n<1,\vspace{7pt}\\
\displaystyle{\frac{N}{dr_n}} & \mbox{if}\quad r_n\ge 1.
\end{array}
\right.
\end{equation}
On the other hand, it follows from \eqref{eq:2.8} and \eqref{eq:4.1} that 
\begin{equation}
\label{eq:4.5}
[S_{K_\theta}(t-s)U^i(s)](x)\le \tilde c_*U^i(x,t),
\quad x\in{\mathbb R}^N,\,\,t>s>0,
\end{equation}
where $\tilde c_*=2c_*C_*$ and $C_*$ is as in \eqref{eq:2.8}.
Combining \eqref{eq:4.3} and \eqref{eq:4.5}, we have
\begin{equation}
\label{eq:4.6}
\begin{split}
 & c_*\int_0^t\int_{{\mathbb R}^N}(t-s)^{-\frac{\ell+j}{d}}K_\theta(x-y,t-s)
 s^{A-\frac{\langle {\bf p}\rangle_n-n}{d}}U^i(y,s)^{|{\bf p}|}\,dy\,ds\\
 & \le c_*(2\gamma_i C)^{|{\bf p}|-1}
\int_0^t(t-s)^{-\frac{\ell+j}{d}} s^{A-\frac{\langle {\bf p}\rangle_n-n}{d}-\kappa(|{\bf p}|-1)}
S_{K_\theta}(t-s)U^i(s)\,ds\\
 & \le c_*\tilde c_*(2\gamma_i C)^{|{\bf p}|-1}U^i(x,t)\int_0^t(t-s)^{-\frac{\ell+j}{d}} 
 s^{A-\frac{\langle {\bf p}\rangle_n-n}{d}-\kappa(|{\bf p}|-1)}\,ds
\end{split}
\end{equation}
for $(x,t)\in{\mathbb R}^N\times(0,T_i)$. 
On the other hand, it follows from (A), \eqref{eq:1.5}, and \eqref{eq:4.4} that 
$$
A-\frac{\langle {\bf p}\rangle_n-n}{d}-\kappa(|{\bf p}|-1)
\ge-1+\frac{\ell+n}{d}>-1.
$$
Since $\ell+j\le\ell+m<d$,
using \eqref{eq:1.5} and \eqref{eq:4.4} again,
we see that 
\begin{equation}
\label{eq:4.7}
\begin{split}
 & \int_0^t(t-s)^{-\frac{\ell+j}{d}} s^{A-\frac{\langle {\bf p}\rangle_n-n}{d}-\kappa(|{\bf p}|-1)}\,ds\\
 & \le Ct^{-\frac{\ell+j}{d}+1+A-\frac{\langle {\bf p}\rangle_n-n}{d}-\kappa(|{\bf p}|-1)}
\le 
\left\{
\begin{array}{ll}
CT_i^{\frac{N(|{\bf p|}-1)}{d}\left(\frac{1}{r_n}-1\right)}t^{-\frac{j-n}{d}} & \mbox{if}\quad r_n<1,\vspace{5pt}\\
Ct^{-\frac{j-n}{d}} & \mbox{if}\quad r_n\ge 1,
\end{array}
\right.\end{split}
\end{equation}
for $0<t<T_i$. 
By \eqref{eq:4.4}, \eqref{eq:4.6}, and \eqref{eq:4.7}, 
taking small enough $\gamma>0$ if necessary, we see that
\begin{equation}
\label{eq:4.8}
\begin{split}
  & c_*\int_0^t\int_{{\mathbb R}^N}(t-s)^{-\frac{\ell+j}{d}}K_\theta(x-y,t-s)s^{A-\frac{\langle {\bf p}\rangle_n-n}{d}}U^i(y,s)^{|{\bf p}|}\,dy\,ds\\
 & \le C\gamma^{|{\bf p}|-1}t^{-\frac{j-n}{d}}U^i(x,t)
 \le\frac{1}{2}t^{-\frac{j-n}{d}}U^i(x,t)
\end{split}
\end{equation}
for $(x,t)\in{\mathbb R}^N\times(0,{T_i})$. 
Therefore, by \eqref{eq:4.2} and \eqref{eq:4.8} 
we apply Lemma~\ref{Lemma:3.1} with $U^\epsilon=U^i$ to find 
$u^i\in C^{m;0}({\mathbb R}^N\times(0,T))$ satisfying \eqref{eq:3.3} with $\epsilon=\epsilon_i$ and
\begin{equation}
\label{eq:4.9}
|\nabla^j u^i(x,t)|\le t^{-\frac{j-n}{d}}U^i(x,t)
\end{equation}
for $(x,t)\in{\mathbb R}^N\times(0,T_i)$ and $j\in\{n,\dots,m\}$.  

Let $\tau\in(0,\min\{T,1\})$. 
Since $T_i>\tau$, 
by \eqref{eq:4.3} and \eqref{eq:4.9}
we see that 
$$
\sup_{i\in\{1,2,\dots\}} \sup_{(x,t)\in{\mathbb R}^N\times[\tau,T_i)}|\nabla^j u^i(x,t)|\le 
\tau^{-\frac{j-n}{d}}\sup_{(x,t)\in{\mathbb R}^N\times[\tau,T_i)}U^i(x,t)\le C<\infty
$$
for $j\in\{n,\dots,m\}$. 
This together with \eqref{eq:3.2} implies that 
$$
\sup_{i\in\{1,2,\dots\}} \sup_{(x,t)\in{\mathbb R}^N\times[\tau,T_i)}|\tilde{F}_{\epsilon_i}(x,t,u^i(x,t),\dots,\nabla^m u^i(x,t))|<\infty.
$$
Applying the parabolic regularity theorems 
(see e.g. \cite{F}*{Chapter~1, Section~3} and \cite{IKK01}*{Section~2}) to integral equation~\eqref{eq:3.3},  
we find $\nu\in(0,1)$ such that 
$$
|\nabla^j u^i(x,t)|\le C,\quad
|\nabla^j u^i(x,t)-\nabla^j u^i(y,s)|\le C(|x-y|^\nu+|t-s|^{\frac{\nu}{d}}),
$$
for $(x,t)$, $(y,s)\in{\mathbb R}^N\times[\tau,T_i)$,  $j\in\{n,\dots,m\}$, and $i\in\{1,2,\dots\}$.
By the Arzel\'a--Ascoli theorem and the diagonal argument 
we find a subsequence $\{u^{i'}\}$ of $\{u^i\}$ 
with $\lim_{i'\to\infty}\epsilon_{i'}=0$ 
and a function $u\in C^{m;0}({\mathbb R}^N\times(0,T))$ such that 
$$
\lim_{i'\to\infty}\sup_{j\in\{n,\dots,m\}}\sup_{(x,t)\in E}|\nabla^j u^{i'}(x,t)-\nabla^j u(x,t)|=0
$$
for compact sets $E$ in ${\mathbb R}^N\times(0,T)$. 
Then, applying the Lebesgue dominated convergence theorem to integral equation~\eqref{eq:3.3}, 
we see that $u$ is a solution to integral equation~(I) 
in ${\mathbb R}^N\times[0,T)$. 
Furthermore, by \eqref{eq:4.3} and \eqref{eq:4.9} 
we observe that $u$ satisfies \eqref{eq:1.7}. 
Thus Theorem~\ref{Theorem:1.2} follows.
$\Box$
\vspace{5pt}

\noindent
{\bf Proof of Corollary~\ref{Corollary:1.1}.}
Corollary~\ref{Corollary:1.1} follows 
Theorem~\ref{Theorem:1.2} with $T=\infty$ and the definition of the Morrey space ${\mathcal M}_{r_n,1}$.
$\Box$
\section{Proofs of Theorems~\ref{Theorem:1.3} and \ref{Theorem:1.4}}
In this section we prove Theorems~\ref{Theorem:1.3} and \ref{Theorem:1.4} by using Lemma~\ref{Lemma:3.1}. 
Furthermore, we prove Corollary~\ref{Corollary:1.2}.
\vspace{5pt}
\newline
{\bf Proof of Theorem~\ref{Theorem:1.4}.}
We can assume, without loss of generality, that $1<q<|\bf p|$. 
Indeed, it follows from the Jensen inequality that
$$
\sup_{x\in{\mathbb R}^N}\biggr(\,\dashint_{B(x,\sigma)}|\phi(y)|^{q'}\,dy\biggr)^{\frac{1}{q'}}
\le\sup_{x\in{\mathbb R}^N}\biggr(\,\dashint_{B(x,\sigma)}|\phi(y)|^q\,dy\biggr)^{\frac{1}{q}}
\quad\mbox{if}\quad 1\le q'\le q.
$$
Let $T\in(0,T_*]$ and assume \eqref{eq:1.9}. 
Let $c_*$ be as in Lemma~\ref{Lemma:3.1}. 
Let $i\in\{1,2,\dots\}$ and fix it. 
Let $T_i:=\min\{T,i\}$ and $\epsilon_i\in(0,1)$ be small enough. 
Then we find $L_i>0$ such that 
$$
c_*L_i\min_{j\in\{0,\dots,m\}}T_i^{-\frac{j}{d}}=\epsilon_i.
$$
Set
\begin{equation}
\label{eq:5.1}
U^i(x,t):=2c_*([S_{K_\theta}(t)\left(|\phi|+L_i\right)^q](x))^\frac{1}{q}.
\end{equation}
Then
\begin{equation}
\label{eq:5.2}
\inf_{(x,t)\in{\mathbb R}^N\times(0,T_i)}t^{-\frac{j}{d}}U^i(x,t)
\ge 2c_* L_i T_i^{-\frac{j}{d}}\ge\epsilon_i,
\quad j\in\{0,\dots,m\}.
\end{equation}
On the other hand, 
it follows from the Jensen inequality, \eqref{eq:1.3}, and \eqref{eq:2.4} that
\begin{equation}
\label{eq:5.3}
c_*[S_{K_\theta}(t)|\phi|](x)
\le c_*([S_{K_\theta}(t)|\phi|^q](x))^\frac{1}{q}
\le\frac{1}{2}U^i(x,t)
\end{equation}
for $(x,t)\in{\mathbb R}^N\times(0,T_i)$.
Since $\lim_{\epsilon_i\to 0}L_i=0$, 
taking small enough $\epsilon_i>0$ if necessary, 
by \eqref{eq:1.9} and \eqref{eq:2.9} we find $C_1>0$ such that
\begin{equation}
\label{eq:5.4}
U^i(x,t)\le C_1\gamma  t^{-\frac{N}{dr_0}}+2c_*L_i
\le 2\gamma  C_1 t^{-\frac{N}{dr_0}},
\quad(x,t)\in{\mathbb R}^N\times(0,T_i). 
\end{equation}
On the other hand, it follows from \eqref{eq:2.8} and \eqref{eq:5.1} that 
\begin{equation}
\label{eq:5.5}
[S_{K_\theta}(t-s)U^i(s)^q](x)\le \tilde c_*U^i(x,t)^q,
\quad x\in{\mathbb R}^N,\,\,t>s>0,
\end{equation}
where $\tilde c_*$ is as in \eqref{eq:4.5}.
Since it follows from \eqref{eq:1.5} that
$$
A-\frac{\langle {\bf p}\rangle_0}{d}-\frac{N(|{\bf p}|-q)}{dr_0}=-1+\frac{N(q-1)}{dr_0}>-1,
$$
by \eqref{eq:5.4} and \eqref{eq:5.5}
we find $C_2>0$ such that
\begin{equation}
\label{eq:5.6}
\begin{split}
 & c_*\int_0^t\int_{{\mathbb R}^N}(t-s)^{-\frac{j}{d}}K_\theta(x-y,t-s)
 s^{A-\frac{\langle {\bf p}\rangle_0}{d}}U^i(y,s)^{|{\bf p}|}\,dy\,ds\\
 & \le c_*(2\gamma C_1)^{|{\bf p}|-q}
\int_0^t(t-s)^{-\frac{j}{d}} s^{A-\frac{\langle {\bf p}\rangle_0}{d}-\frac{N(|{\bf p}|-q)}{dr_0}}
\int_{{\mathbb R}^N}K_\theta(x-y,t-s)U^i(y,s)^q\,dy\,ds\\
 & \le c_*\tilde c_*M(2\gamma C_1)^{|{\bf p}|-q}U^i(x,t)^q\int_0^t
 (t-s)^{-\frac{j}{d}} s^{A-\frac{\langle {\bf p}\rangle_0}{d}-\frac{N(|{\bf p}|-q)}{dr_0}}\,ds\\
 & \le c_*\tilde c_*(2\gamma C_1)^{|{\bf p}|-1}t^{-\frac{N(q-1)}{dr_0}}U^i(x,t)
 \int_0^t(t-s)^{-\frac{j}{d}} s^{A-\frac{\langle {\bf p}\rangle_0}{d}-\frac{N(|{\bf p}|-q)}{dr_0}}\,ds\\
 & \le C_2\gamma^{|{\bf p}|-1}t^{-\frac{j}{d}}U^i(x,t)
\end{split}
\end{equation}
for $(x,t)\in{\mathbb R}^N\times(0,T_i)$. 
Then, taking small enough $\gamma>0$ if necessary, 
we have 
\begin{equation}
\label{eq:5.7}
 c_*\int_0^t\int_{{\mathbb R}^N}(t-s)^{-\frac{j}{d}}K_\theta(x-y,t-s)s^{A-\frac{\langle {\bf p}\rangle_0}{d}}U^i(y,s)^{|{\bf p}|}\,dy\,ds
\le\frac{1}{2}t^{-\frac{j}{d}}U^i(x,t)
\end{equation}
for $(x,t)\in{\mathbb R}^N\times(0,T_i)$. 
By \eqref{eq:5.2}, \eqref{eq:5.3}, and \eqref{eq:5.7}  
we apply Lemma~\ref{Lemma:3.1} with $U^\epsilon=U^i$ 
to find $u^i\in C^{m;0}({\mathbb R}^N\times(0,T_i))$
satisfying \eqref{eq:3.3} with $\epsilon=\epsilon_i$ and
\begin{equation}
\label{eq:5.8}
|\nabla^j u^i(x,t)|\le t^{-\frac{j}{d}}U^i(x,t),
\quad j\in\{0,\dots,m\},
\end{equation}
for $(x,t)\in{\mathbb R}^N\times(0,T_i)$.  
Then, by the same arguments as in the proof of Theorem~\ref{Theorem:1.2} 
we find a solution~$u\in C^{m;0}({\mathbb R}^N\times(0,T))$ to integral equation~(I) 
in ${\mathbb R}^N\times[0,T)$. 
Furthermore, 
by \eqref{eq:5.4} and \eqref{eq:5.8}
we see that $u$ satisfies \eqref{eq:1.10}. 
Thus Theorem~\ref{Theorem:1.4} follows.
$\Box$
\vspace{5pt}

\noindent
{\bf Proof of Corollary~\ref{Corollary:1.2}.}
Similarly to the proof of Corollary~\ref{Corollary:1.1}, 
Corollary~\ref{Corollary:1.2} follows 
Theorem~\ref{Theorem:1.4} with $T=\infty$ 
and the definition of the Morrey space ${\mathcal M}_{r_0,q}$.
$\Box$
\vspace{5pt}

\noindent{\bf Proof of Theorem~\ref{Theorem:1.3}.}
Let $c_*$ be as in Lemma~\ref{Lemma:3.1}. 
It follows from \eqref{eq:1.8} that 
$$
\sup_{x\in{\mathbb R}^N}\int_{B(x,t^{\frac{1}{d}})}|\phi(y)|\,dy
\le \sup_{x\in{\mathbb R}^N}\int_{B(x,T^{\frac{1}{d}})}|\phi(y)|\,dy
\le\gamma T^{\frac{N}{d}-\frac{N}{dr_0}},
\quad 0<t\le T.
$$
This together with \eqref{eq:2.9} implies that 
\begin{equation}
\label{eq:5.9}
\|S_{K_\theta}(t)|\phi|\|_{L^\infty}\le C\gamma T^{\frac{N}{d}-\frac{N}{dr_0}}t^{-\frac{N}{d}},
\quad t\in(0,T).
\end{equation}
Let $\epsilon>0$. Then we find $L_\epsilon>0$ such that 
$$
c_*L_\epsilon\min_{j\in\{0,\dots,m\}}T^{-\frac{j}{d}}=\epsilon. 
$$
Set 
\begin{equation}
\label{eq:5.10}
U^\epsilon(x,t):=2c_*[S_{K_\theta}(t)(|\phi|+L_\epsilon)](x).
\end{equation}
These imply that
\begin{equation*}
\begin{split}
 & \inf_{(x,t)\in{\mathbb R}^N\times(0,T]}t^{-\frac{j}{d}}U^\epsilon(x,t)
\ge\epsilon,
\quad j\in\{0,\dots,m\},\\
 & c_*[S_{K_\theta}(t)|\phi|](x)
\le\frac{1}{2}U^\epsilon(x,t),
\quad (x,t)\in{\mathbb R}^N\times(0,T).
\end{split}
\end{equation*}
Furthermore, taking small enough $\epsilon>0$ if necessary, 
by \eqref{eq:5.9} we find $C_1>0$ such that
\begin{equation}
\label{eq:5.11}
U^\epsilon(x,t)\le C_1\gamma T^{\frac{N}{d}-\frac{N}{dr_0}} t^{-\frac{N}{d}}+2c_*L_\epsilon
\le 2C_1 \gamma  T^{\frac{N}{d}-\frac{N}{dr_0}} t^{-\frac{N}{d}},
\quad
(x,t)\in{\mathbb R}^N\times(0,T).
\end{equation}
On the other hand, it follows from \eqref{eq:2.8} and \eqref{eq:5.10} that 
\begin{equation}
\label{eq:5.12}
[S_{K_\theta}(t-s)U^\epsilon(s)](x)\le \tilde c_*U^\epsilon(x,t),
\quad x\in{\mathbb R}^N,\,\,t>s>0,
\end{equation}
where $\tilde c_*$ is as in \eqref{eq:4.5}.
Since it follows from \eqref{eq:1.5} and (B) that
$$
A-\frac{\langle {\bf p}\rangle_0}{d}-\frac{N(|\bf p|-1)}{d}
=-1+\frac{N(|\bf p|-1)}{d}\bigg(\frac{1}{r_0}-1\bigg)>-1,
$$
similarly to \eqref{eq:5.6},
by \eqref{eq:5.11} and \eqref{eq:5.12} we find $C_2>0$ such that
\begin{equation*}
\begin{split}
 & c_*\int_0^t\int_{{\mathbb R}^N}(t-s)^{-\frac{j}{d}}K_\theta(x-y,t-s)s^{A-\frac{\langle {\bf p}\rangle_0}{d}}U^\epsilon(y,s)^{|{\bf p}|}\,dy\,ds\\
 & \le c_*(2C_1\gamma T^{\frac{N}{d}-\frac{N}{dr_0}})^{|{\bf p}|-1}
 \int_0^t(t-s)^{-\frac{j}{d}} s^{A-\frac{\langle {\bf p}\rangle_0}{d}-\frac{N(|{\bf p}|-1)}{d}}
\int_{{\mathbb R}^N}K_\theta(x-y,t-s)U^\epsilon(y,s)\,dy\,ds\\
 & \le C_2\gamma^{|{\bf p}|-1}t^{-\frac{j}{d}}U^\epsilon(x,t)
\end{split}
\end{equation*}
for $(x,t)\in{\mathbb R}^N\times(0,T)$ and $j\in\{0,\dots,m\}$.
Then, applying the same arguments as in the proof of Theorem~\ref{Theorem:1.4}
with $q$ replaced by $1$, 
we complete the proof of Theorem~\ref{Theorem:1.3}.
$\Box$
%
\section{Proof of Theorem~\ref{Theorem:1.5}}
We prove Theorem~\ref{Theorem:1.5} and Corollary~\ref{Corollary:1.3}. 
\vspace{5pt}
\newline
{\bf Proof of Theorem~\ref{Theorem:1.5}.}
Let $M\ge e$ and set $\Phi_M(s):=s[\log(M+s)]^\beta$ for $s>0$. 
Then, taking large enough $M\ge e$ if necessary, 
we have:
\begin{itemize}
  \item[{\rm (i)}] 
  $\Phi_M$ is convex in $(0,\infty)$;
  \item[{\rm (ii)}] 
  The function $(0,\infty)\ni s\mapsto$
  $s^{\frac{|\bf p|-1}{2}}[\log(M+s)]^{-\beta |\bf p|}$ is monotone increasing.
\end{itemize}
Furthermore, by \eqref{eq:1.11} we see that
\begin{equation}
\label{eq:6.1}
\begin{split}
C^{-1}\Phi_M(s) & \le\Phi(s)\le C\Phi_M(s),\\
C^{-1}s[\log(M+s)]^{-\beta} & \le\Phi_M^{-1}(s)\le Cs[\log(M+s)]^{-\beta},
\end{split}
\end{equation}
for $s>0$. 
It follows from \eqref{eq:1.12} and \eqref{eq:6.1} that 
\begin{equation}
\label{eq:6.2}
\sup_{x\in\mathbb R^N}\Phi_M^{-1}\left[\,\dashint_{B(x,\sigma)}
\Phi_M(T^{\frac{N}{d}}|\phi(y)|)\,dy\,\right]\le C\gamma \rho(\sigma T^{-\frac{1}{d}}),
\quad 0<\sigma<T^{\frac{1}{d}}.
\end{equation}
Set 
\begin{equation}
\label{eq:6.3}
V(x,t):=[S_{K_\theta}(t)\Phi_M(T^{\frac{N}{d}}|\phi|)](x),
\quad
\tau(t):=T^{-\frac{1}{d}}t^{\frac{1}{d}}.
\end{equation}
Then, by \eqref{eq:2.9} and \eqref{eq:6.2} we see that
\begin{equation}
\label{eq:6.4}
\begin{split}
{\|V(t)\|_{L^\infty}}
 & \le Ct^{-\frac{N}{d}}\sup_{x\in{\mathbb R}^N}\int_{B(x,t^{1/d})}
\Phi_M(T^{\frac{N}{d}}|\phi(y)|)\,dy\\
 & \le C\Phi_M(C\gamma \rho(\tau(t)))
 \le C\gamma \rho(\tau(t))[\log(M+C\gamma \rho(\tau(t)))]^\beta\\
 & \le C\gamma \tau(t)^{-N}\left|\log\frac{\tau(t)}{2}\right|^{-\frac{N}{d(1+A)-\langle{\bf p}\rangle_0}+\beta}=:\gamma  \xi(\tau(t))
\end{split}
\end{equation}
for $t\in(0,T)$. 
Here the last inequality in \eqref{eq:6.4} follows from \eqref{eq:1.11} and 
\begin{equation*}
\begin{split}
\rho(\tau)[\log(M+C\rho(\tau))]^\beta
 & =O\left(\tau^{-N}|\log \tau|^{-\frac{N}{d(1+A)-\langle{\bf p}\rangle_0}}|\log \tau|^\beta\right)\\
 & =O\left(\tau^{-N}|\log \tau|^{-\frac{N}{d(1+A)-\langle{\bf p}\rangle_0}+\beta}\right)
\quad\mbox{as}\quad\tau\to +0.
\end{split}
\end{equation*}

Let $c_*$ be as in Lemma~\ref{Lemma:3.1}. 
For any $\epsilon>0$, 
let $L_\epsilon>0$ be such that 
$$
c_*L_\epsilon\min_{j\in\{0,\dots,m\}}T^{-\frac{N+j}{d}}=\epsilon. 
$$
Then, taking small enough $\epsilon>0$ if necessary, 
by \eqref{eq:6.4} we see that
\begin{equation}
\label{eq:6.5}
\|V(t)\|_{L^\infty}+\Phi_M(L_\epsilon)\le 2\gamma \xi(\tau(t))\quad\mbox{for}\quad t\in (0,T).
\end{equation}
Set
\begin{equation}
\label{eq:6.6}
V^\epsilon(x,t):=V(x,t)+\Phi_M(L_\epsilon),\quad 
U^\epsilon(x,t):=2c_*T^{-\frac{N}{d}}\Phi_M^{-1}\left(V^\epsilon(x,t)\right).
\end{equation}
Then 
\begin{equation}
\label{eq:6.7}
\inf_{(x,t)\in{\mathbb R}^N\times(0,T)}t^{-\frac{j}{d}}U^\epsilon(x,t)\ge 2c_*T^{-\frac{N+j}{d}}L_\epsilon\ge\epsilon.
\end{equation}
On the other hand, 
by \eqref{eq:1.3}, \eqref{eq:2.4}, and \eqref{eq:6.6} 
we apply the Jensen inequality to obtain 
\begin{equation}
\label{eq:6.8}
c_*[S_{K_\theta}(t)|\phi|](x)
\le c_*T^{-\frac{N}{d}}\Phi_M^{-1}\left[S_{K_\theta}(t)\Phi_M\left(T^{\frac{N}{d}}|\phi|\right)\right](x)
\le\frac{1}{2}U^\epsilon(x,t)
\end{equation}
for $(x,t)\in{\mathbb R}^N\times(0,T)$.

We can assume, without loss of generality, that $\gamma \le 1/2$. 
By property~(ii), \eqref{eq:6.1}, \eqref{eq:6.4}, \eqref{eq:6.5}, and \eqref{eq:6.6} we have
\begin{align*}
0 \le \frac{U^\epsilon(x,t)^{|{\bf p}|}}{V^\epsilon(x,t)}
&
=(2c_*)^{|\bf p|}T^{-\frac{N|{\bf p}|}{d}}\frac{[\Phi_M^{-1}(V^\epsilon(x,t))]^{|{\bf p}|}}{V^\epsilon(x,t)}\\
 & \le CT^{-\frac{N|{\bf p}|}{d}}V^\epsilon(x,t)^{|{\bf p}|-1}[\log(M+V^\epsilon(x,t))]^{-\beta |{\bf p}|}\\
 & =CT^{-\frac{N|{\bf p}|}{d}}V^\epsilon(x,t)^{\frac{|{\bf p}|-1}{2}}V^\epsilon(x,t)^{\frac{|{\bf p}|-1}{2}}[\log(M+V^\epsilon(x,t))]^{-\beta |{\bf p}|}\\
 & \le CT^{-\frac{N|{\bf p}|}{d}}(2\gamma\xi(\tau(t)))^{\frac{|{\bf p}|-1}{2}}
 (2\gamma \xi(\tau(t)))^{\frac{|{\bf p}|-1}{2}}
 [\log(M+2\gamma \xi(\tau(t)))]^{-\beta |{\bf p}|}
 \\
 & \le C\gamma^{\frac{|{\bf p}|-1}{2}}
 T^{-\frac{N|{\bf p}|}{d}}\xi(\tau(t))^{|{\bf p}|-1}[\log(M+\xi(\tau(t)))]^{-\beta |{\bf p}|}
\end{align*}
for $(x,t)\in{\mathbb R}^N\times(0,T)$. 
This together with (D) and the definition of $\xi(\tau(t))$ implies that 
\begin{equation}
\label{eq:6.9}
\begin{split}
0\le \frac{U^\epsilon(x,t)^{|{\bf p}|}}{V^\epsilon(x,t)}
 & \le C\gamma^{\frac{|{\bf p}|-1}{2}} T^{-\frac{N|{\bf p}|}{d}}\tau(t)^{-N(|{\bf p}|-1)}
 \left|\log\frac{\tau(t)}{2}\right|^{-\frac{N(|{\bf p}|-1)}{d(1+A)-\langle {\bf p}\rangle_0}
 +\beta(|{\bf p}|-1)-\beta |{\bf p}|}\\
 & =C\gamma^{\frac{|{\bf p}|-1}{2}} T^{-\frac{N|{\bf p}|}{d}}\tau(t)^{-N(|{\bf p}|-1)}\left|\log\frac{\tau(t)}{2}\right|^{-1-\beta}
\end{split}
\end{equation}
for $(x,t)\in{\mathbb R}^N\times(0,T)$. 
Similarly, 
\begin{equation}
\label{eq:6.10}
\begin{split}
\frac{V^\epsilon(x,t)}{U^\epsilon(x,t)} & =(2c_*)^{-1}T^{\frac{N}{d}}
\frac{V^\epsilon(x,t)}{\Phi_M^{-1}(V^\epsilon(x,t))}\\
 & \le CT^{\frac{N}{d}}[\log(M+V^\epsilon(x,t))]^\beta
\le CT^{\frac{N}{d}}\biggr|\log\frac{\tau(t)}{2}\biggr|^\beta
\end{split}
\end{equation}
for $(x,t)\in{\mathbb R}^N\times(0,T)$. 
On the other hand, it follows from \eqref{eq:2.8}, \eqref{eq:6.3}, and \eqref{eq:6.6} that 
\begin{equation}
\label{eq:6.11}
[S_{K_\theta}(t-s)V^\epsilon(s)](x)\le C_*V^\epsilon(x,t),
\quad x\in{\mathbb R}^N,\,\,t>s>0,
\end{equation}
where $C_*$ is as in \eqref{eq:2.8}.
Combining \eqref{eq:6.9}, \eqref{eq:6.10}, and \eqref{eq:6.11} 
we obtain 
\begin{equation}
\label{eq:6.12}
\begin{split}
 & c_*\int_0^t\int_{{\mathbb R}^N}(t-s)^{-\frac{j}{d}}K_\theta(x-y,t-s)
 s^{A-\frac{\langle {\bf p}\rangle_0}{d}}U^\epsilon(y,s)^{|{\bf p}|}\,dy\,ds\\
 & \le C\gamma^{\frac{|{\bf p}|-1}{2}}T^{-\frac{N|{\bf p}|}{d}}
  \\
 & \qquad
 \times
 \int_0^t(t-s)^{-\frac{j}{d}}s^{A-\frac{\langle {\bf p}\rangle_0}{d}}
 \tau(s)^{-N(|{\bf p}|-1)}\left|\log\frac{\tau(s)}{2}\right|^{-1-\beta}
 S_{K_\theta}(t-s)V^\epsilon(s)\,ds\\
 & \le CC_*\gamma^{\frac{|{\bf p}|-1}{2}}T^{-\frac{N|{\bf p}|}{d}}V^\epsilon(x,t)
 \int_0^t(t-s)^{-\frac{j}{d}}s^{A-\frac{\langle {\bf p}\rangle_0}{d}}
 \tau(s)^{-N(|{\bf p}|-1)}\left|\log\frac{\tau(s)}{2}\right|^{-1-\beta}\,ds\\
 & \le C\gamma^{\frac{|{\bf p}|-1}{2}}T^{-\frac{N(|{\bf p}|-1)}{d}}U^\epsilon(x,t)\\
 & \qquad\qquad
 \times\left|\log\frac{\tau(t)}{2}\right|^{\beta}
 \int_0^t(t-s)^{-\frac{j}{d}}s^{A-\frac{\langle {\bf p}\rangle_0}{d}}
 \tau(s)^{-N(|{\bf p}|-1)}\left|\log\frac{\tau(s)}{2}\right|^{-1-\beta}\,ds
\end{split}
\end{equation} 
for $(x,t)\in{\mathbb R}^N\times(0,T)$. 
Since it follows from (D) and \eqref{eq:6.3} that
\begin{equation*}
\begin{split}
 & \int_0^t(t-s)^{-\frac{j}{d}}s^{A-\frac{\langle {\bf p}\rangle_0}{d}}\tau(s)^{-N(|{\bf p}|-1)}\left|\log\frac{\tau(s)}{2}\right|^{-1-\beta}\,ds\\
 & =T^{\frac{N(|{\bf p}|-1)}{d}}\int_0^t(t-s)^{-\frac{j}{d}}s^{-1}\left|\log\frac{\tau(s)}{2}\right|^{-1-\beta}\,ds
 \le CT^{\frac{N(|{\bf p}|-1)}{d}}t^{-\frac{j}{d}}\left|\log\frac{\tau(t)}{2}\right|^{-\beta}
\end{split}
\end{equation*}
for $t\in(0,T)$, taking small enough $\gamma>0$ if necessary, 
we deduce from \eqref{eq:6.12}  
that
\begin{equation}
\label{eq:6.13}
\begin{split}
 & c_*\int_0^t\int_{{\mathbb R}^N}(t-s)^{-\frac{j}{d}}K_\theta(x-y,t-s)s^{-\frac{\langle {\bf p}\rangle_0}{d}}U^\epsilon(y,s)^{|{\bf p}|}\,dy\,ds\\
 & \le C\gamma^{\frac{|{\bf p}|-1}{2}}t^{-\frac{j}{d}}U^\epsilon(x,t)
 \le\frac{1}{2}t^{-\frac{j}{d}}U^\epsilon(x,t)
\end{split}
\end{equation} 
for $(x,t)\in\mathbb R^N\times(0,T)$ and $j=\{0,\dots,m\}$.
Therefore, by \eqref{eq:6.7}, \eqref{eq:6.8}, and \eqref{eq:6.13} 
we apply Lemma~\ref{Lemma:3.1} to find $u^\epsilon\in C^{m;0}({\mathbb R}^N\times(0,T))$
satisfying \eqref{eq:3.3} and
\begin{equation}
\label{eq:6.14}
|\nabla^j u^\epsilon(x,t)|\le t^{-\frac{j}{d}}U^\epsilon(x,t),\quad j\in\{0,\dots,m\},
\end{equation} 
for $(x,t)\in{\mathbb R}^N\times(0,T)$.  
Then, by the same arguments as in the proof of Theorem~\ref{Theorem:1.2} 
we find a solution~$u\in C^{m;0}({\mathbb R}^N\times(0,T))$ to integral equation~(I) in ${\mathbb R}^N\times[0,T)$. 
Furthermore, 
by \eqref{eq:6.1}, \eqref{eq:6.3}, and \eqref{eq:6.14}
we have 
\begin{equation*}
\begin{split}
|\nabla^j u^\epsilon(x,t)| & \le CT^{-\frac{N}{d}}t^{-\frac{j-n}{d}}\Phi_M^{-1}(V^\epsilon(x,t))\\
 & \le CT^{-\frac{N}{d}}t^{-\frac{j}{d}}\xi(\tau(t))[\log(M+\xi(\tau(t)))]^{-\beta}
\le Ct^{-\frac{N+j}{d}}\biggr|\log\frac{t}{2T}\biggr|^{-\frac{N}{d}}
\end{split}
\end{equation*}
for $(x,t)\in{\mathbb R}^N\times(0,T)$. 
This implies inequality~\eqref{eq:1.13}.
Thus Theorem~\ref{Theorem:1.4} follows. 
$\Box$
\vspace{5pt}

\noindent
{\bf Proof of Corollary~\ref{Corollary:1.3}.}
Let 
$$
\beta=\frac{N}{d(1+A)-\langle{\bf p}\rangle_0}.
$$
Then
$$
\Phi(\gamma \rho(\sigma T^{-\frac{1}{d}}))=\gamma \rho(\sigma T^{-\frac{1}{d}})
\left[\log(e+\gamma \rho(\sigma T^{-\frac{1}{d}}))\right]^{\frac{N}{d(1+A)-\langle{\bf p}\rangle_0}}
\le C\gamma T^{\frac{N}{d}}\sigma^{-N}
$$
for $0<\sigma<T^{\frac{1}{d}}$. 
Then Corollary~\ref{Corollary:1.3}. follows from Theorem~\ref{Theorem:1.4}.
$\Box$
\section{Applications}
It is known that 
fundamental solutions to a large class of linear parabolic operators satisfy condition~(G) 
for some $T_*\in(0,\infty]$, $d>0$, and $L>0$,
and our main results are applicable to the Cauchy problem
for various nonlinear parabolic equations. 
In this section we focus on the Cauchy problem for 
nonlinear parabolic equation~\eqref{eq:1.1} with 
$$
{\mathcal L}=(-\Delta)^{\frac{d}{2}},\qquad d>0, 
$$
and show the validity and the advantage of our main results. 
\subsection{Semilinear parabolic equations}
Consider the Cauchy problem for a semilinear parabolic equation
\begin{equation}
\tag{SP}
\left\{
\begin{array}{ll}
\partial_t u+(-\Delta)^{\frac{d}{2}}u=|u|^p,\quad & x\in{\mathbb R}^N,\,\,t>0,\vspace{3pt}\\
u(x,0)=\phi(x), & x\in{\mathbb R}^N,
\end{array}
\right.
\end{equation}
where $d>0$ and $p>1$.
Then 
$$
\ell=n=m=A=0,\quad |{\bf p}|=p,\quad \langle {\bf p}\rangle_0=0.
$$ 
Applying Theorems~\ref{Theorem:1.3}, \ref{Theorem:1.4}, and \ref{Theorem:1.5} to problem~(SP), 
we have: 
\begin{theorem}
\label{Theorem:7.1}
Consider Cauchy problem~{\rm (SP)}, where $d>0$ and $p>1$. 
Then the same statements as in 
Theorems~{\rm\ref{Theorem:1.3}}, {\rm\ref{Theorem:1.4}}, and {\rm\ref{Theorem:1.5}} hold 
in the cases $r_0<1$, $r_0>1$, and $r_0=1$, respectively, where
$$
r_0=\frac{N(p-1)}{d}.
$$
\end{theorem}
In the case either $0<d\le 2$ or $d\in\{4,6,\dots\}$, 
Theorem~\ref{Theorem:7.1} has been already proved in \cite{HI01} and \cite{IKO}, respectively, 
and it is shown that sufficient conditions in 
Theorems~\ref{Theorem:1.3}, \ref{Theorem:1.4}, and~\ref{Theorem:1.5} are sharp. 
See also Section~1. 
\subsection{Viscous Hamilton-Jacobi equations}
Consider the Cauchy problem for a viscous Hamilton-Jacobi equation
\begin{equation}
\tag{VHJ}
\left\{
\begin{array}{ll}
\partial_t u+(-\Delta)^{\frac{d}{2}}u=|\nabla u|^p,\quad & x\in{\mathbb R}^N,\,\,t>0,\vspace{3pt}\\
u(x,0)=\phi(x), & x\in{\mathbb R}^N,
\end{array}
\right.
\end{equation}
where $d>1$ and $p>1$.
Then
\begin{equation*}
\begin{split}
 & n\in\{0,1\},\quad \ell=A=0,\quad m=1,\quad |{\bf p}|=p,
\quad \langle {\bf p}\rangle_0=p,\quad \langle {\bf p}\rangle_1=1,\\
 & r_0=\frac{N(p-1)}{d-p},\quad r_1=\frac{N(p-1)}{d-1}.
\end{split}
\end{equation*}
It follows that $r_0<1$ if and only if $p<p_{HJ}$, 
where
$$
p_{HJ}:=\frac{N+d}{N+1}\in(1,d).
$$
Applying Theorems~\ref{Theorem:1.2}, \ref{Theorem:1.3}, \ref{Theorem:1.4}, 
and \ref{Theorem:1.5} to problem~(VHJ), 
we have: 
\begin{theorem}
\label{Theorem:7.2}
Let $p>1$ and $d>1$. 
\begin{itemize}
  \item[{\rm (a)}] 
  Let $1<p<p_{HJ}$. Then problem~{\rm (VHJ)} possesses a local-in-time solution if 
  $\phi\in L^1_{{\rm uloc}}$.
  \item[{\rm (b)}] 
  Let $p=p_{HJ}$. Then the same statement as in Theorem~{\rm\ref{Theorem:1.5}} holds with 
  $$
  \rho(s)=s^{-N}[\log(e+s^{-1})]^{-\frac{N}{d-p}}.
  $$
  In particular, there exists $C_1>0$ such that, if  
  $$
  |\phi(x)|\le C_1|x|^{-N}\biggr[\log\biggr(e+\frac{1}{|x|}\biggr)\biggr]^{-\frac{N}{d}-1}+C,\quad x\in{\mathbb R}^N,
  $$
  for some $C>0$, 
  problem~{\rm (VHJ)} possesses a local-in-time solution.
  \item[{\rm (c)}] 
  Let $p_{HJ}<p<d$ and $q>1$. Then there exists $C_2>0$ such that, if 
  $$
  \sup_{x\in{\mathbb R}^N}\,\sup_{0<\sigma<T^{\frac{1}{d}}}\,
  \sigma^{\frac{d-p}{p-1}}\,
  \biggr(\,\dashint_{B(x,\sigma)}|\phi(y)|^q\,dy\biggr)^{\frac{1}{q}}\le C_2
  $$
  for some $T\in(0,\infty]$, 
  problem~{\rm (VHJ)} possesses a solution in ${\mathbb R}^N\times[0,T)$.
  \item[{\rm (d)}] Let $p>1$. 
  Then there exists $C_3>0$ such that, if  
  $$
  \sup_{x\in{\mathbb R}^N}\,\sup_{0<\sigma<T^{\frac{1}{d}}}\,
   \sigma^{\frac{d-1}{p-1}}\,\dashint_{B(x,\sigma)}|\nabla\phi(y)|\,dy\le C_3
  $$
  for some $T\in(0,\infty]$, 
  problem~{\rm (VHJ)} possesses a solution in ${\mathbb R}^N\times[0,T)$.
\end{itemize}
\end{theorem}
We remark that, 
in the case of $1<d\le 2$, 
thanks to \cite{AB}, \cite{DI}, and \cite{KW},
Definition~\ref{Definition:1.1} implies that 
if problem~(VHJ) possesses a local-in-time solution~$u$, then 
the solution~$u$ can be extended as a global-in-time solution to (VHJ).
For the case of $d=2$, 
due to \cite{BSW},  
the well-posedness of local-in-time solutions holds
in the following cases:
\begin{itemize}
  \item $\phi\in L^q$ for $q>r_1\ge1$ or $q=r_1>1$;
  \item $\phi\in L^1$ if $p<p_{HJ}$;
  \item $\phi\in W^{1,q}(\mathbb R^N)$ if $p\in[1,\infty)$ and $q>r_1\ge1$ or $q=r_1>1$.
\end{itemize}
(See also \cites{A, BL}.) 
Comparing with these results,
we see that Theorem~\ref{Theorem:7.2} includes some new criterions for the existence of solutions to problem~(VHJ) 
even in the case of $d=2$, in particular, in assertions~(b) and (c). 
\subsection{Nonlinear parabolic equations with $\ell>0$ and $m=0$}
Consider the Cauchy problem
\begin{equation}
\tag{gCD}
\left\{
\begin{array}{ll}
\partial_t u+(-\Delta)^{\frac{d}{2}}u=
\displaystyle{\sum_{|\alpha|=\ell}}a_\alpha\partial_x^\alpha (|u|^{p-1}u),\quad & x\in{\mathbb R}^N,\,\,t>0,\vspace{3pt}\\
u(x,0)=\phi(x), & x\in{\mathbb R}^N,
\end{array}
\right.
\end{equation}
where $\ell\in\{1,2,\dots\}$, $0<\ell<d$, $\{a_\alpha\}\subset{\mathbb R}$, and $p>1$.
Then 
$$
n=A=0,\quad |{\bf p}|=p,\quad\langle {\bf p}\rangle_0=0,
\quad
r_0=\frac{N(p-1)}{d-\ell}.
$$
Applying Theorem~\ref{Theorem:1.2} to problem~(gCD), 
we have: 
\begin{theorem}
\label{Theorem:7.3}
Let $\ell\in\{1,2,\dots\}$, $0<\ell<d$, $\{a_\alpha\}\subset{\mathbb R}$, and $p>1$. 
Then there exists $\gamma>0$ such that, if
$$
\sup_{x\in{\mathbb R}^N}\,\sup_{0<\sigma<T^{\frac{1}{d}}}\,
\sigma^{\frac{d-\ell}{p-1}}\,\dashint_{B(x,\sigma)}|\phi(y)|\,dy\le\gamma
$$
for some $T\in(0,\infty]$, then problem~{\rm (gCD)} possesses a solution in ${\mathbb R}^N\times[0,T)$.
\end{theorem}
In the case when $1<d\le 2$ and $\ell=1$, 
the comparison principle holds for problem~(gCD).
Then, similarly to problem~(VHJ), 
Definition~\ref{Definition:1.1} implies that 
if problem~(gCD) possesses a local-in-time solution~$u$, then 
the solution~$u$ can be extended as a global-in-time solution to (gCD). 

Problem~(gCD) is a generalization of the Cauchy problem for a convection-diffusion equation
\begin{equation}
\tag{CD}
\left\{
\begin{array}{ll}
\partial_t u+(-\Delta)^{\frac{d}{2}} u={\bf a}\cdot\nabla(|u|^{p-1}u),\quad & x\in{\mathbb R}^N,\,\,t>0,\vspace{3pt}\\
u(x,0)=\phi(x), & x\in{\mathbb R}^N,
\end{array}
\right.
\end{equation}
where $0<d\le 2$, $p>1$, and ${\bf a}\in{\mathbb R}^N$. 
The solvability and the asymptotic behavior of solutions to problem~(CD) 
have been studied in many papers (see e.g. \cites{Carpio, EZ01, EZ02, IKK02, IwaK, HOS} and references therein). 
Theorem~\ref{Theorem:7.3} also includes some new criterions for the existence of solutions even to problem~(CD). 
(Compare with \cite{HOS}.)
%
\subsection{Nonlinear parabolic equations with $\ell>0$ and $m=1$}
Consider the Cauchy problem for a higher-order parabolic equation with gradient nonlinearity 
\begin{equation}
\tag{HG}
\left\{
\begin{array}{ll}
\partial_t u+(-\Delta)^{\frac{d}{2}}u
=-\nabla\cdot(|\nabla u|^{p-1}\nabla u),\quad & x\in{\mathbb R}^N,\,\,t>0,\vspace{3pt}\\
u(x,0)=\phi(x), & x\in{\mathbb R}^N,
\end{array}
\right.
\end{equation}
where $d>2$ and $p>1$.
Then 
\begin{equation*}
\begin{split}
 & n\in\{0,1\},\quad A=0,\quad \ell=1,\quad |{\bf p}|=p,\quad 
 \langle {\bf p}\rangle_0=p,\quad \langle {\bf p}\rangle_1=1,\\
 & r_0=\frac{N(p-1)}{d-p-1},\quad r_1=\frac{N(p-1)}{d-2}.
\end{split}
\end{equation*}
Applying Theorem~\ref{Theorem:1.2} to problem~(HG), we have:
\begin{theorem}
\label{Theorem:7.4}
Let $d>2$ and $p>1$. 
\begin{itemize}
  \item[{\rm (a)}] 
  Let $1<p\le d-1$. 
  Then there exists $C_1>0$ such that, if
  $$
  \sup_{x\in{\mathbb R}^N}\,\sup_{0<\sigma<T^{\frac{1}{d}}}\,
  \sigma^{\frac{d-p-1}{p-1}}\,
  \dashint_{B(x,\sigma)}|\phi(y)|\,dy\le C_1
  $$
  for some $T\in(0,\infty]$, then problem~{\rm (HG)} possesses a solution in ${\mathbb R}^N\times[0,T)$.
  \item[{\rm (b)}] 
  Let $p>1$. 
  Then there exists $C_2>0$ such that, if
  $$
  \sup_{x\in{\mathbb R}^N}\,\sup_{0<\sigma<T^{\frac{1}{d}}}\,
  \sigma^{\frac{d-2}{p-1}}\,\dashint_{B(x,\sigma)}|\nabla \phi(y)|\,dy\le C_2
  $$
  for some $T\in(0,\infty]$, then problem~{\rm (HG)} possesses a solution in ${\mathbb R}^N\times[0,T)$.
\end{itemize}
\end{theorem}
Problem~(HG) with $d=4$ appears 
in the study of mathematical models describing epitaxial growth of thin film 
(see e.g. \cites{EGP, EGK, KSW, IMO, ORS} and references therein for related results). 
In \cite{IMO} the authors gave sufficient conditions 
for the existence of local-in-time solutions and global-in-time solutions 
by the use of uniformly local weak Lebesgue spaces. 
Theorem~\ref{Theorem:7.4} with $d=4$ improves their results.
\medskip

\noindent
{\bf Acknowledgment.} 
The authors of this paper were supported in part 
by JSPS KAKENHI Grant Number JP19H05599.
The second author and the third author were supported in part by JSPS KAKENHI Grant Numbers JP20K03689 
and JP20KK0057, respectively. 
  
\end{document}